\pdfoutput=1
\RequirePackage{ifpdf}
\ifpdf 
\documentclass[pdftex]{sigma}
\else
\documentclass{sigma}
\fi

\numberwithin{equation}{section}
\newtheorem{Theorem}{Theorem}[section]
\newtheorem{Lemma}[Theorem]{Lemma}
\newtheorem{Proposition}[Theorem]{Proposition}
{ \theoremstyle{definition}
\newtheorem{Definition}[Theorem]{Definition}
\newtheorem{Notation}[Theorem]{Notation}
\newtheorem{Remark}[Theorem]{Remark} }
\DeclareMathOperator{\SP}{sp}

\begin{document}
\newcommand{\nablad}{\nabla\hskip -1.05ex \vrule width.1ex height 1.3ex depth -.6ex \hskip 1.05ex}

\newcommand{\arXivNumber}{math.OA/0602212}

\allowdisplaybreaks

\renewcommand{\thefootnote}{$\star$}

\renewcommand{\PaperNumber}{082}

\FirstPageHeading

\ShortArticleName{Locally Compact Quantum Groups.
A~von Neumann Algebra Approach}

\ArticleName{Locally Compact Quantum Groups.\\
A~von Neumann Algebra Approach\footnote{This paper is a~contribution to the Special Issue on Noncommutative Geometry and
Quantum Groups in honor of Marc A.~Rief\/fel.
The full collection is available at
\href{http://www.emis.de/journals/SIGMA/Rieffel.html}{http://www.emis.de/journals/SIGMA/Rieffel.html}}}

\Author{Alfons VAN DAELE}

\AuthorNameForHeading{A.~Van Daele}

\Address{Department of Mathematics, University of Leuven, Celestijnenlaan 200B,\\
B-3001 Heverlee, Belgium}
\Email{\href{mailto:alfons.vandaele@wis.kuleuven.be}{alfons.vandaele@wis.kuleuven.be}}

\ArticleDates{Received February 06, 2014, in f\/inal form July 28, 2014; Published online August 05, 2014}

\Abstract{In this paper, we give an alternative approach to the theory of locally compact quantum groups, as developed
by Kustermans and Vaes.
We start with a~von Neumann algebra and a~comultiplication on this von Neumann algebra.
We assume that there exist faithful left and right Haar weights.
Then we develop the theory within this von Neumann algebra setting.
In~[\textit{Math. Scand.} \textbf{92} (2003), 68--92] locally compact quantum groups are also studied in the von Neumann algebraic context.
This approach is independent of the original $C^*$-algebraic approach in the sense that the earlier results are not
used.
However, this paper
is not really independent because for many proofs, the reader is referred to the original
paper where the $C^*$-version is developed.
In this paper, we give a~completely self-contained approach.
Moreover, at various points, we do things dif\/ferently.
We have a~dif\/ferent treatment of the antipode.
It is similar to the original treatment in~[\textit{Ann. Sci.
  \'Ecole Norm. Sup.~(4)} \textbf{33} (2000), 837--934].
But together with the fact that we work in the von Neumann algebra framework, it allows us to use an idea
from~[\textit{Rev. Roumaine Math. Pures Appl.} \textbf{21} (1976),
  1411--1449] to obtain the uniqueness of the Haar weights in an early stage.
We take advantage of this fact when deriving the other main results in the theory.
We also give a~slightly dif\/ferent approach to duality.
Finally, we collect, in a~systematic way, several important formulas.
In an appendix, we indicate very brief\/ly how the $C^*$-approach and the von Neumann algebra approach eventually yield
the same objects.
The passage from the von Neumann algebra setting to the $C^*$-algebra setting is more or less standard.
For the other direction, we use a~new method.
It is based on the observation that the Haar weights on the $C^*$-algebra extend to weights on the double dual with
central support and that all these supports are the same.
Of course, we get the von Neumann algebra by cutting down the double dual with this unique support projection in the
center.
All together, we see that there are many advantages when we develop the theory of locally compact quantum groups in the
von Neumann algebra framework, rather than in the $C^*$-algebra framework.
It is not only simpler, the theory of weights on von Neumann algebras is better known and one needs very little to go
from the $C^*$-algebras to the von Neumann algebras.
Moreover, in many cases when constructing examples, the von Neumann algebra with the coproduct is constructed from the
very beginning and the Haar weights are constructed as weights on this von Neumann algebra (using left Hilbert algebra
theory).
This paper is written in a~concise way.
In many cases, only indications for the proofs of the results are given.
This information should be enough to see that these results are correct.
We will give more details in forthcoming paper,
which will be expository, aimed at non-specialists.
See also~[\textit{Bull. Kerala Math. Assoc.}  (2005),
  153--177] for an `expanded' version of the appendix.}


\Classification{26L10; 16L05; 43A99}

\renewcommand{\thefootnote}{\arabic{footnote}}
\setcounter{footnote}{0}

\section{Introduction}

Let~$M$ be a~von Neumann algebra and~$\Delta$ a~comultiplication on~$M$ (see Def\/inition~\ref{1.1}
for
a~precise def\/inition).
The pair $(M,\Delta)$ is called a~{\it locally compact quantum group} (in the von Neumann algebraic sense) if there
exist faithful left and right Haar weights (see Def\/inition~\ref{2.1}).
This def\/inition is due to Kustermans and Vaes (see~\cite{K-V3}).

In their fundamental papers~\cite{K-V1} and~\cite{K-V2}, Kustermans and Vaes develop the theory of locally compact
quantum groups in the $C^*$-algebraic framework and in~\cite{K-V3}, they show that both the original $C^*$-algebra
approach and the von Neumann algebra approach give the same objects.
There is indeed a~standard procedure to go from a~locally compact quantum group in the $C^*$-algebra setting to
a~locally compact quantum group in the von Neumann algebra sense (and vice versa).

In this paper, we present an alternative approach to the theory of locally compact quantum groups.
The basic dif\/ference is that we develop the main theory in the framework of von Neumann algebras (and not in the
$C^*$-algebraic setting as was done in~\cite{K-V1} and~\cite{K-V2}).
It is well-known that the von Neumann algebra setting is, in general, simpler to work in.
To begin with, the def\/inition of a~locally compact quantum group in the von Neumann algebra framework is already less
complicated than in the $C^*$-algebra setting.
Also the theory of weights on von Neumann algebras is better known than the theory of weights on $C^*$-algebras.
One has to be a~bit more careful with using the various topologies, but on the other hand, one does not have to worry
about multiplier algebras.

We also describe a~(relatively) quick way to go from a~locally compact quantum group in the $C^*$-algebraic sense to
a~locally compact quantum group in the von Neumann algebraic sense.
We do not need to develop the $C^*$-theory to do this.
In some sense, this justif\/ies our choice.
Remark that the other direction, from von Neumann algebras to $C^*$-algebras, is the easier one (and standard).

However, the dif\/ference of this work with the other (earlier) approaches not only lies in the fact that we work in the
von Neumann algebra framework.
We also have a~slightly dif\/ferent approach to construct the antipode (see Section~\ref{Section2}).
We do not use operator space techniques (as e.g.\ in~\cite{V-VD}).
This in combination with the use of Connes' cocycle Radon--Nikodym theorem (an idea that we found in earlier works by
Stratila, Voiculescu and Zsido, see~\cite{S-V-Z1,S-V-Z2,S-V-Z3}) allows us to obtain uniqueness of
the Haar weights in an earlier stage of the development.
This in turn will yield other simplif\/ications.

Finally, our approach in this paper is also self-contained.
In their paper on the von Neumann algebraic approach~\cite{K-V3}, Kustermans and Vaes do not really use results from the
earlier paper on the $C^*$-algebra approach~\cite{K-V2}, but nevertheless, it is hard to read it without the f\/irst paper
because of the fundamental references to this f\/irst paper.
In fact, we also rely less on results from other papers (e.g.\ on weights on $C^*$-algebras or about manageability of
multiplicative unitaries) as is done in the original works.
Also for the proofs that are omitted, this is the case.
Moreover, where possible, we avoid working with unbounded operators and weights (more than in the original papers) but
we try to use bounded operators and normal linear functionals.
We do not use operator valued weights at all.

The {\it content of the paper} is as follows.
In Section~\ref{Section2} we work with a~von Neumann algebra and a~comultiplication.
We consider the antipode~$S$, together with an involutive operator~$K$ on the Hilbert space that implements the antipode
in the sense that, roughly speaking, $S(x)^*=KxK$ when $x\in \mathcal D(S)$.
The right Haar weight is needed to construct this operator~$K$ and the left Haar weight is used to prove that it is
densely def\/ined.
This last property is closely related with the right regular representation being unitary.
Also in this section, we focus on various other densities.
It makes this section longer than the others, but the reason for doing so is that these density results are closely
related with the construction of the antipode and the operator~$K$.
Finally, in this section, we modify the def\/inition of the antipode so that it becomes more tractable.
We discuss its polar decomposition (with the scaling group $(\tau_t)$ and the unitary antipode~$R$) and we prove the
basic formulas about this modif\/ied antipode, needed further in the paper.

Our approach here is not so very dif\/ferent from the way this is done by Kustermans and Vaes in the sense that we use the same ideas.
Among other things, we also use Kustermans' trick to prove that the right regular representation is unitary.
A~sound knowledge of various aspects of the Tomita--Takesaki theory and its relation with weights on von Neumann algebras
is necessary for understanding the arguments.
However, we avoid the use of operator valued weighs.
We will include the necessary background in the notes we plan to write~\cite{VD8} but including more details here
would make this paper too long.

In Section~\ref{Section3} we give the main results.
One of these results is the uniqueness of the Haar weights.
The formulas involving the scaling group $(\tau_t)$ and the unitary antipode~$R$, proven in Section~\ref{Section2}, together with
Connes' cocycle Radon--Nikodym theorem, are used to show the uniqueness.
Here, our approach is quite dif\/ferent from the original one in~\cite{K-V1} and~\cite{K-V2} and uses an idea found in~\cite{S-V-Z2}.
From the uniqueness, and again using basic formulas involving the scaling group and the unitary antipode from Section~\ref{Section2},
the main results are relatively easy to prove.

In Section~\ref{Section4} we treat the dual.
This is more or less standard.
Our approach is again slightly dif\/ferent in the way we use the results obtained in Sections~\ref{Section2} and~\ref{Section3}.
Also, since we are basically only considering the von Neumann algebra version, the construction of the dual is somewhat simpler.

In Section~\ref{Section5} we collect a~set of formulas.
The main ingredients are the various objects (the left and right Haar weights with their modular structures, the left
and right regular representations, the antipode with the scaling group and the unitary antipode, the operators on the
Hilbert space implementing these automorphism groups, \dots), for the original pair $(M,\Delta)$, as well as for the dual
pair $(\widehat M,\widehat\Delta)$.
In fact, this section and these formulas can well be used as a~fairly complete chapter needed to work with locally compact quantum groups.

In Section~\ref{Section6} we draw some conclusions and discuss possible further research along the lines of this paper.

We have chosen to discuss the procedure to pass from the $C^*$-algebraic locally compact quantum groups to the von
Neumann algebraic ones in Appendix~\ref{appendixA}.
We do this because it is not really needed for the development of the theory as it is done in this paper.
Here again, our approach is rather dif\/ferent from the original one.
The main idea is to pass f\/irst to the double dual $A^{**}$ of the $C^*$-algebra~$A$.
Then it is quickly proven that the supports of the invariant weights (of the type used in this theory), are all the same
central projection in the double dual.
Cutting down the double dual by this central projection gives us the von Neumann algebra~$M$.
The coproduct, as well as the Haar weights, are obtained by f\/irst extending the corresponding original objects to the
double dual and then restricting them again to this von Neumann algebra~$M$.
The converse is standard and of course makes use of the results in the paper.
A~more or less independent treatment of the connection of the two approaches is found in an expanded version of this
appendix, see~\cite{VD7}.
We also use Appendix~\ref{appendixA} to say something more about the relation of our approach with the one by Masuda, Nakagami and
Woronowicz (in~\cite{M-N} and~\cite{M-N-W}).

In this paper, we will not give full details.
We give precise def\/initions and statements, but often we will only {\it sketch proofs}.
We give suf\/f\/iciently many details so that the reader should be `convinced' about the result.
In fact any reader, familiar with the Tomita--Takesaki theory in relation with the theory of weights on von Neumann
algebras, should be able to complete the proofs without too much ef\/fort.
On the other hand, for the less experienced reader, we refer to a~forthcoming paper {\it Notes on locally compact
quantum groups}~\cite{VD8}.
These notes are intended as lecture notes on the subject, for (young) researchers who want to learn about locally
compact quantum groups.
So full details of the proofs of the results in this paper will be found there.
This style of writing allows us to keep this paper reasonable in size, while on the other hand, we still are able to
make it, to a~great extent, self-contained.
In Appendix~\ref{appendixA} we give even less details because this is not so essential for the development here.
In~\cite{VD7} an expanded version of this appendix is found, but again, for more details we refer to~\cite{VD8}.
Finally let us also refer to the book of Timmermann~\cite{Tim} for an overview of the theory of multiplier Hopf
$*$-algebras and algebraic quantum groups, within the context of the theory of locally compact quantum groups.
Much information about the purely algebraic theory can be found there and this can be helpful to understand the
technically far more dif\/f\/icult analytical theory.
Also the original papers on the theory of multiplier Hopf ($*$)-algebras~\cite{VD2} and~\cite{VD4}, although certainly
not necessary for understanding this paper, can be helpful.

We would like to emphasize the {\it importance of the original work} by Kustermans and Vaes, also for this alternative treatment.
We do not really use results from their work, but certainly we have been greatly inspired by their results and techniques.
Without their pioneering work, this paper would not have been written.
It is worthwhile mentioning that the PhD Thesis of Vaes~\cite{V2} (for those who have access to this work) is easier to
read than the original paper~\cite{K-V2}.
Also the paper by Masuda, Nakagami and Woronowicz~\cite{M-N-W}, treating independently the theory of locally compact
quantum groups, has helped us to develop our new approach.
Throughout the paper, we will not always repeat to refer to the original works, but we will do so where we feel this is appropriate.

Let us now f\/inish this introduction with some {\it basic references} and {\it standard notations} used in this paper.
We will also say something about the dif\/ference in conventions used in the f\/ield.

When $\mathcal H$ is a~Hilbert space, we will use $\mathcal B(\mathcal H)$ to denote the von Neumann algebra of all
bounded linear operators on $\mathcal H$.
We use $M_*$ for the space of normal linear functionals on a~von Neumann algebra~$M$ and in particular $\mathcal
B(\mathcal H)_*$ for normal linear functionals on $\mathcal B(\mathcal H)$.
When~$\omega$ is a~such a~functional, we use $\overline\omega$ for the linear functional def\/ined~by
$\overline\omega(x)=\omega(x^*)^-=\overline{\omega(x^*)}$ (where in all these cases, the $^-$ stands for complex conjugation).
On one occasion, we will also need the absolute value $|\omega|$ and the norm $\| \omega \|$ of a~normal linear functional.
If $\xi$ and~$\eta$ are two vectors in the Hilbert space $\mathcal H$, we will write $\langle\,\cdot\,\xi,\eta\rangle$
to denote the normal linear functional~$\omega$ on $\mathcal B(\mathcal H)$ given by $x\mapsto \langle x\xi,\eta\rangle$.
In this case we have e.g.\ $\overline\omega(x)=\langle x\eta,\xi\rangle$ and, provided $\|\xi\|=\|\eta\|=1$, that
$|\omega|(x)=\langle x\xi,\xi\rangle$ and $\|\omega\|=1$.

We refer to~\cite{S-Z,T2,T3} for the theory of $C^*$-algebras and von Neumann algebras.

We will work with normal semi-f\/inite weights on von Neumann algebras and with lower semi-continuous densely def\/ined
weights on $C^*$-algebras.
We will use the standard notations for the objects associated with such weights.
If e.g.~$\psi$ is a~normal semi-f\/inite weight on a~von Neumann algebra~$M$, we will use $\mathcal N_\psi$ for the left
ideal of elements $x\in M$ such that $\psi(x^*x)<\infty$.
Also $\mathcal M_\psi$ will be the hereditary $*$-subalgebra $\mathcal N_\psi^* \mathcal N_\psi$ of~$M$, spanned by the
elements $x^*y$ with $x,y\in \mathcal N_\psi$.
We will use the G.N.S.-representation associated with such a~weight.
The Hilbert space will be denoted by $\mathcal H_\psi$ while $\Lambda_\psi$ is used for the canonical map from $\mathcal
N_\psi$ to the space $\mathcal H_\psi$.
We will let the von Neumann algebra act directly on its G.N.S.
space, i.e.\ we will drop the notation $\pi_\psi$.
The modular operator on $\mathcal H_\psi$ (in the case of a~faithful weight) will be denoted by $\nabla$ (and {\it not}
by~$\Delta$ because we reserve~$\Delta$ for comultiplications).
We will use $(\sigma_t^\psi)$ for the modular automorphisms.

Again we refer to~\cite{T3} for the theory of weights on $C^*$-algebras and von Neumann algebras, as well as for the
modular theory and its relation with weights.
See also~\cite{St}.
For the original work on left Hilbert algebras, there is of course~\cite{T1}.

We will be using various (continuous) one-parameter groups of automorphisms.
We assume~$\sigma$-weak continuity, but one can easily see that~$\sigma$-weak continuity for one-parameter groups of
$*$-automorphisms implies also continuity for the stronger operator topologies (like the~$\sigma$-strong or even
the~$\sigma$-strong-$*$ topology).
We also have that the map $t\mapsto\omega\circ\alpha_t$ is continuous for any $\omega\in M_*$ with the norm topology on
$M_*$ when~$\alpha$ is a~continuous one-parameter group of automorphisms.
We will be interested in analytical elements and the analytical generator $\alpha_{\frac{i}{2}}$ of such a~one-parameter
group~$\alpha$.
A~few things can be found in Chapter~VIII of~\cite{T3} and a~nice reference is also Appendix~F in~\cite{M-N-W}.
It seems to be better to f\/irst def\/ine the analytical generator for the action of $\mathbb R$ on $M_*$, dual to the
one-parameter group (because this is norm continuous), and def\/ine the analytical generator on~$M$ by taking the adjoint.
Doing so, we get that the linear map $\alpha_{\frac{i}{2}}$ is closed for the~$\sigma$-weak topology and that the
analytical elements form a~core with respect to the~$\sigma$-strong-$*$ topology.

When we write the tensor product of spaces, we will always mean completed tensor products.
In the case of two Hilbert spaces, this is the Hilbert space tensor product.
In the case of $C^*$-algebras, it is understood to be the minimal $C^*$-tensor product.
Finally, for von Neumann algebras we take the usual von Neumann algebra (i.e.\ the spatial) tensor product.

Unfortunately, there are a~number of dif\/ferent conventions used in this f\/ield by dif\/ferent authors/schools.
In Hopf algebras e.g.\ it is common to endow the dual of a~(f\/inite-dimensional) Hopf algebra with a~coproduct simply~by
dualizing the product (see e.g.~\cite{A}) whereas in the theory of locally compact quantum groups, usually the opposite
coproduct on the dual is taken.
In the earlier works on Kac algebras (see~\cite{E-S}), the left regular representation is def\/ined as the adjoint of what
is commonly used now.
Kustermans and Vaes work mainly with the left regular representation (as in the case of Kac algebras), whereas Baaj and
Skandalis (in~\cite{B-S}) and Masuda, Nakagami and Woronowicz (in~\cite{M-N-W}) prefer the right regular representation
as their starting point.
Also a~dif\/ferent convention in~\cite{M-N-W} is used for the polar decomposition of the antipode.

In this paper, we will mainly follow the conventions used by Kustermans and Vaes in their original papers.
In a~few occasions, mostly as a~consequence of the dif\/ference in approach, we will choose slightly dif\/ferent conventions.
In that case, we will clearly say so.
It will only be the case in the process of obtaining the main results.
In the formulation of the main results, we will be in accordance with the conventions in the papers of Kustermans and Vaes.
Also the papers~\cite{VD5} and~\cite{K-V4} are interesting as short survey papers.

\section{The antipode: construction and properties. Densities}\label{Section2}

Let~$M$ be a~von Neumann algebra and denote by $M\otimes M$ the von Neumann tensor product of~$M$ with itself.
Recall the following def\/inition which is the basic ingredient of this paper.
Let us also assume that~$M$ acts on the Hilbert space $\mathcal H$ in standard form.

\begin{Definition}\label{1.1}
Let~$\Delta$ be a~unital and normal $*$-homomorphism from~$M$ to $M\otimes M$.
Then~$\Delta$ is called a~{\it comultiplication} on~$M$ if $(\Delta\otimes\iota)\Delta = (\iota\otimes\Delta)\Delta$
(coassociativity), where~$\iota$ is used to denote the identity map from~$M$ to itself.
\end{Definition}

The standard example comes from a~locally compact group~$G$.
We take $M=L^\infty(G)$ and def\/ine~$\Delta$ on~$M$ by $\Delta(f)(r,s)=f(rs)$ whenever $f\in L^\infty(G)$ and $r,s\in G$.
We identify $M\otimes M$ with $L^\infty(G\times G)$.

\medskip

{\bf A~preliminary def\/inition of the antipode and f\/irst properties.}
We f\/irst def\/ine the following subspace of the von Neumann algebra~$M$.
\begin{Definition}\label{1.2}
For an element $x\in M$ we say that $x\in\mathcal D_0$ if there is an element $x_1 \in M$ satisfying the following condition:
\begin{itemize}\itemsep=0pt
\item For all $\varepsilon > 0$ and vectors $\xi_1, \xi_2, \dots, \xi_n,
\eta_1, \eta_2, \dots, \eta_n$ in $\mathcal H$,
there exist elements
$p_1, p_2, \dots$, $p_m,q_1, q_2, \dots, q_m$ in~$M$ such that
\begin{gather*}
 \Big\|x \xi_k\otimes \eta_k - \textstyle\sum\limits_j \Delta(p_j)(\xi_k\otimes q_j^*\eta_k)\Big\| < \varepsilon,
\qquad
\Big\|x_1 \xi_k\otimes \eta_k - \textstyle\sum\limits_j \Delta(q_j)(\xi_k\otimes p_j^*\eta_k)\Big\| < \varepsilon
\end{gather*}
 for all~$k$.
\end{itemize}
\end{Definition}

We will see later that, because of forthcoming assumptions, we will have $x_1=0$ if $x=0$.
Therefore it will be possible to def\/ine a~linear map $S_0$ on $\mathcal D_0$ by letting $S_0(x)=x_1^*$.
Then for this map one can prove the following properties:
\begin{enumerate}\itemsep=0pt
\item[i)] if $x\in \mathcal D_0$, then $S_0(x)^*\in \mathcal D_0$ and
$S_0(S_0(x)^*)^*=x$,

\item[ii)] if $x,y\in \mathcal D_0$, then $xy\in \mathcal D_0$ and $S_0(xy)=S_0(y)S_0(x)$,

\item[iii)] the map $x\to S_0(x)^*$ is closed for the strong operator topology on~$M$.
\end{enumerate}

The properties i) and iii) are immediate consequences of the def\/inition while~ii) is obtained using a~simple calculation.
We refer to similar arguments given in the proofs of Propositions~\ref{1.9} and~\ref{1.10}.
However, remark that we will not really use these results for~$S_0$.

This operator would be a~candidate for the antipode (see the following remarks), but we will not def\/ine the antipode
like this but rather through its polar decomposition (see the Def\/initions~\ref{1.22} and~\ref{1.23} later in this section).
It is expected that the two def\/initions coincide, but we have not been able to show this.
Fortunately, it is not necessary for the further development in this paper.

\begin{Remark}
\quad
\begin{enumerate}\itemsep=0pt
\item[i)] The def\/inition of $S_0$ above is inspired by a~result in Hopf ($*$-)algebra theory (see~\cite{A} and~\cite{Sw}).
Indeed if $(H,\Delta)$ is a~Hopf algebra with antipode~$S$ and if $a\in H$, then using the Sweedler notation, we get
\begin{gather*}
a\otimes 1  = \textstyle\sum\limits_{(a)} a_{(1)} \otimes a_{(2)}S(a_{(3)})
 = \textstyle\sum\limits_{(a)} \Delta(a_{(1)})(1\otimes S(a_{(2)})).
\end{gather*}
So, if we have a~Hopf $*$-algebra and if we write $\sum\limits_j p_j\otimes q_j = \sum\limits_{(a)} a_{(1)}\otimes
S(a_{(2)})^*$, we get
\begin{gather*}
\textstyle\sum\limits_j \Delta(q_j)(1\otimes p_j^*)  = \textstyle\sum\limits_{(a)} \Delta(S(a_{(2)})^*)(1\otimes
a_{(1)}^*)
 = \textstyle\sum\limits_{(a)} S(a_{(3)})^* \otimes S(a_{(2)})^*a_{(1)}^*
 = S(a)^*\otimes 1.
\end{gather*}
If we do not have a~$*$-structure, we have a~similar formula, but we loose the symmetry.

\item[ii)] This formula can also be illustrated in the case $M=L^\infty(G)$ where~$G$ is a~locally compact group.
In this case we know that $S(f)(r)=f(r^{-1})$ when $f\in L^\infty(G)$ and $r\in G$.
If we approximate
\begin{gather*}
f(r)=f\big(rs\cdot s^{-1}\big)\simeq \textstyle\sum\limits_i p_i(rs)\overline{q_i(s^{-1})},
\end{gather*}
we get
\begin{gather*}
\textstyle\sum\limits_i q_i(rs) \overline{p_i(s)}\simeq \overline{f(s\cdot(rs)^{-1})}=\overline{f(r^{-1})}.
\end{gather*}

\item[iii)] We have used the above idea in the construction of the antipode for Hopf $C^*$-algebras in~\cite{V-VD}.
In fact, also the construction of the antipode in the paper~\cite{K-V2} uses this idea, but that is less obvious.

\item[iv)] The well-def\/inedness of $S_0$ in the general case is a~problem and it is also not clear whether or not there are
even non-trivial elements in~$\mathcal D_0$.
As we will see later in this section, the left and right Haar weights will be used to solve this problem.

\item[v)] One of the nice aspects however of this approach to the antipode is that {\it it does not depend on the possible
choices} of the left and the right Haar weights.
\end{enumerate}
\end{Remark}

The reader should have these remarks in mind further in this section.

\medskip

{\bf The involutive operator~$\boldsymbol{K}$ implementing the antipode.}
In what follows, we will see how the existence of the Haar weights eventually leads to, not only the {\it well-definedness}
of this preliminary antipode $S_0$, but also gives the {\it density of the domain} $\mathcal D_0$.
But as we already mentioned, we will not def\/ine the antipode in this way.
On the other hand, we will do something similar and def\/ine a~map like $S_0(\,\cdot\,)^*$, but on the Hilbert space level.

To do this, we will {\it now assume the existence of a~right Haar weight}.
We recall the def\/inition (see e.g.~\cite{K-V3}).

\begin{Definition}\label{1.4}
Let~$M$ be a~von Neumann algebra and~$\Delta$ a~comultiplication on~$M$ (as in Def\/inition~\ref{1.1}).
A~{\it right Haar weight} on~$M$ is a~faithful, normal semi-f\/inite weight on~$M$ such that
\begin{gather*}
\psi((\iota\otimes \omega)\Delta(x))=\omega(1)\psi(x),
\end{gather*}
whenever $x\in M$, $x\geq 0$ and $\psi(x)<\infty$ and when $\omega \in M_*$ and $\omega \geq 0$ (right invariance).
\end{Definition}

We will now further in this section {\it fix a~right Haar weight}~$\psi$.

We consider the G.N.S.-representation of~$M$ for~$\psi$.
Let $\mathcal N_\psi$ be the set of elements $x\in M$ such that $\psi(x^*x)<\infty$.
We will use $\Lambda_\psi$ to denote the canonical map from the $\mathcal N_\psi$ to $\mathcal H_\psi$.
As usual, we extend~$\psi$ to the $*$-subalgebra $\mathcal M_\psi$ (def\/ined as $\mathcal N_\psi^*\mathcal N_\psi)$.
We have $\langle \Lambda_\psi(x)$, $\Lambda_\psi(y)\rangle=\psi(y^*x)$ for all $x,y\in \mathcal N_\psi$.
We consider~$M$ as acting directly on $\mathcal H_\psi$ (i.e.\ we drop the notation $\pi_\psi$) and so we will write
$x\Lambda_\psi(y)=\Lambda_\psi(xy)$ when $x\in M$ and $y\in\mathcal N_\psi$.

We refer to~\cite{St} and~\cite{T3} for details about weights and the G.N.S.-construction for weights.

Also by now, the construction of the {\it right regular representation} has become standard.
We recall it here and refer to e.g.~\cite{K-V2} and~\cite{K-V3} as well as to~\cite{E-S} (and also~\cite{VD8}) for
details.

\begin{Proposition}\label{1.5}
There exists a~bounded operator~$V$ from $\mathcal H_\psi \otimes \mathcal H$ to itself, characterized $($and defined$)$~by
\begin{gather*}
((\iota\otimes\omega)V)\Lambda_\psi(x)=\Lambda_\psi((\iota\otimes\omega)\Delta(x)),
\end{gather*}
whenever $x\in \mathcal N_\psi$ and $\omega\in \mathcal B(\mathcal H)_*$.
It has its `second leg' in~$M$, i.e.\ $V\in \mathcal B(\mathcal H_\psi)\otimes M$ and satisfies the following formulas:
\begin{enumerate}\itemsep=0pt
\item[$i)$] $V^*V=1$ $($i.e.~$V$ is an isometry$)$,

\item[$ii)$] $V(x\otimes 1)= \Delta(x)V$ for all $x\in M$,

\item[$iii)$] $(\iota\otimes \Delta)V=V_{12}V_{13}$ $($where we use the standard `leg numbering' notation$)$.
\end{enumerate}
\end{Proposition}

Roughly speaking we have $V(\Lambda_\psi(x)\otimes\xi)=\sum\limits_{(x)}\Lambda_\psi(x_{(1)})\otimes x_{(2)}\xi$ when we
use the Sweedler notation $\Delta(x)=\sum\limits_{(x)} x_{(1)}\otimes x_{(2)}$ for $x\in M$.

Recall that the invariance is used to get that $(\iota\otimes\omega)\Delta(x)\in \mathcal N_\psi$ when $x\in \mathcal N_\psi$,
a~result which is needed to def\/ine~$V$ as above.
Also it is the right invariance that implies that~$V$ is an isometry.
It is known that in general it seems impossible to show that~$V$ is a~unitary without further assumptions.
In our approach, we will get unitarity in some sense as a~`byproduct' of the further study of the antipode (see
Proposition~\ref{1.15}).
Because we do not yet know that~$V$ is unitary, we need to formulate condition ii) as we have done and we can not (yet)
write $\Delta(x)=V(x\otimes 1)V^*$.
This will follow later.

It is easy to show that in the case $M=L^\infty(G)$, the operator~$V$ is indeed intimately related with the right
regular representation of~$G$ on $L^2(G)$ (with the right Haar measure on~$G$).
Recall that the right Haar weight on $L^\infty(G)$ is obtained by integration with respect to the right Haar measure on~$G$.

The {\it next step} is to construct the operator $S_0(\,\cdot\,)^*$ on the Hilbert space level, in this case, on~$\mathcal H_\psi$.
It will be denoted by~$K$.
The def\/inition is very much as in Def\/inition~\ref{1.2}.

\begin{Definition}\label{1.6}
Let $\xi\in \mathcal H_\psi$.
We say that $\xi\in \mathcal D(K)$ if there is a~vector $\xi_1\in \mathcal H_\psi$ satisfying the following condition:

For all $\varepsilon > 0$ and vectors $\eta_1, \eta_2, \dots, \eta_n$ in $\mathcal H$, there exist elements
$p_1, p_2, \dots, p_m$, $q_1, q_2, \dots, q_m$ in $\mathcal N_\psi$ such that
\begin{gather*}
 \Big\|\xi\otimes \eta_k - V\Big(\textstyle\sum\limits_j \Lambda_\psi(p_j) \otimes q_j^*\eta_k\Big)\Big\| < \varepsilon,
\qquad
 \Big\| \xi_1\otimes \eta_k - V\Big(\textstyle\sum\limits_j \Lambda_\psi(q_j)\otimes p_j^*\eta_k\Big)\Big\| < \varepsilon
\end{gather*}
for all~$k$.
\end{Definition}

Remark that this def\/inition is indeed similar to Def\/inition~\ref{1.2} because, roughly speaking, the operator~$V$ is the
map $p\otimes q^*\mapsto \Delta(p)(1\otimes q^*)$ on the Hilbert space level.

Again, we would like to def\/ine the operator~$K$ by $K\xi=\xi_1$ but we need the following result.

\begin{Lemma}
Let $\xi$ and $\xi_1$ be as Definition~{\rm \ref{1.6}} and assume $\xi=0$.
Then also $\xi_1=0$.
\end{Lemma}

\begin{proof}
In this proof, we will take for $\mathcal H$ the space $\mathcal H_\psi$ with the G.N.S.-representation of~$M$.

Take vectors $\eta_1$ and $\eta_2$ in $\mathcal H_\psi$ and $\varepsilon>0$.
By assumption we have elements $(p_j)$ and $(q_j)$ in~$\mathcal N_\psi$ so that
\begin{gather}
 \Big\| \textstyle\sum\limits_j \Lambda_\psi(p_j) \otimes q_j^*\eta_1\Big\| < \varepsilon,
\qquad
\Big \| \xi_1\otimes \eta_2 - V\Big(\textstyle\sum\limits_j \Lambda_\psi(q_j)\otimes p_j^*\eta_2\Big)\Big\| < \varepsilon.
\label{(1.2)}
\end{gather}
Recall that by assumption $\xi=0$ and that~$V$ is isometric.

Now take any pair $\rho_1$, $\rho_2$ of right bounded vectors in $\mathcal H_\psi$.
Recall that a~vector $\rho\in \mathcal H_\psi$ is called right bounded if there is a~bounded operator, necessarily
unique and denoted as $\pi'(\rho)$, satisfying $x\rho=\pi'(\rho)\Lambda_\psi(x)$ for all $x\in \mathcal N_\psi$.

Then we have
\begin{gather*}
\textstyle\sum\limits_j \langle \Lambda_\psi(p_j)\otimes q_j^*\eta_1, \pi'(\rho_1)^*\eta_2\otimes \rho_2 \rangle
 =\textstyle\sum\limits_j \langle \pi'(\rho_1)\Lambda_\psi(p_j)\otimes \eta_1, \eta_2\otimes q_j\rho_2 \rangle
\\
\qquad{}
 =\textstyle\sum\limits_j \langle p_j\rho_1\otimes \eta_1, \eta_2\otimes \pi'(\rho_2)\Lambda_\psi(q_j)\rangle
 =\textstyle\sum\limits_j \langle \rho_1\otimes \pi'(\rho_2)^*\eta_1, p_j^*\eta_2\otimes \Lambda_\psi(q_j)\rangle.
\end{gather*}
It follows that
\begin{gather*}
\Big|\textstyle\sum\limits_j \langle \Lambda_\psi(q_j)\otimes p_j^*\eta_2,\pi'(\rho_2)^*\eta_1\otimes \rho_1 \rangle \Big|  \leq
\Big\|\textstyle\sum\limits_j \Lambda_\psi(p_j)\otimes q_j^*\eta_1\Big\| \|\pi'(\rho_1)^*\eta_1 \otimes \rho_2\|
\\
\qquad{}
 \leq \varepsilon \| \pi'(\rho_1)^*\eta_1\| \|\rho_2\|.
\end{gather*}
This implies that
\begin{gather*}
| \langle \xi_1\otimes \eta_2, V(\pi'(\rho_2)^*\eta_1\otimes \rho_1)\rangle |  \leq \Big\| \xi_1\otimes \eta_2 -
V\Big(\textstyle\sum\limits_j \Lambda_\psi(q_j)\otimes p_j^*\eta_2\Big) \Big\| \| V(\pi'(\rho_2)^* \eta_1 \otimes \rho_1) \|
\\
\qquad
\phantom{\leq}{}
+ \Big| \Big\langle V\Big(\textstyle\sum\limits_j   \Lambda_\psi(q_j)\otimes p_j^*\eta_2\Big),V(\pi'(\rho_2)^* \eta_1 \otimes
\rho_1)\Big\rangle \Big|
\\
\qquad{}
 \leq \varepsilon \| \pi'(\rho_2)^* \eta_1 \| \|\rho_1 \| + \varepsilon \| \pi'(\rho_1)^*\eta_1\| \| \rho_2\|.
\end{gather*}
This is true for all~$\varepsilon$.
Therefore we have
\begin{gather*}
\langle \xi_1\otimes \eta_2, V(\pi'(\rho_2)^*\eta_1\otimes \rho_1)\rangle =0
\end{gather*}
for all right bounded vectors $\rho_1$, $\rho_2$ and all $\eta_1$ in $\mathcal H_\psi$.
Because the set of vectors
$\pi'(\rho_2)^*\eta_1\otimes \rho_1$ span a~dense subspace of $\mathcal H_\psi\otimes \mathcal H_\psi$, we see that
$\xi_1\otimes \eta_2$ is orthogonal to the range of~$V$.
But as it clearly also belongs to the range of~$V$ (as will follow from~\eqref{(1.2)} above), it has to be zero.
Hence $\xi_1=0$.
This completes the proof.
\end{proof}

This argument is not fundamentally dif\/ferent from a~similar argument in~\cite{K-V2}.

\begin{Definition}\label{1.8}
If $\xi\in \mathcal D(K)$ and if $\xi_1$ is as in Def\/inition~\ref{1.6}, we set $K\xi=\xi_1$.
\end{Definition}

Remark that our operator~$K$ is essentially the operator $G^*$ in the work of Kustermans and Vaes (as we will see later~-- cf.\
e.g.\ Remark~\ref{4.10}).
Therefore it should not be a~surprise that the techniques used above to def\/ine~$K$ and to show that it is well-def\/ined
are similar as those used in~\cite{K-V2}.
Observe that we will not use the symbol~$G$ for this operator as this is commonly used to denote a~locally compact group.

Just as in the case of Def\/inition~\ref{1.2}, we get easily the following results.

\begin{Proposition}\label{1.9}\quad
\begin{enumerate}\itemsep=0pt
\item[$i)$] If $\xi\in\mathcal D(K)$, then $K\xi\in\mathcal D(K)$ and $K(K\xi)=\xi$.

\item[$ii)$] $K$ is a~closed operator.
\end{enumerate}
\end{Proposition}

\begin{proof}
i) This is immediately clear from the symmetry we have in Def\/inition~\ref{1.6}.

ii) Assume that we have a~sequence $(\xi_i)$ in $\mathcal D(K)$ and two vectors $\xi$, $\xi'$ in $\mathcal H_\psi$ so that
$\xi_i\to \xi$ and $K\xi_i\to \xi'$.
We have to show that $\xi\in \mathcal D(K)$ and $K\xi=\xi'$.
In other words, we must verify that the pair $(\xi,\xi')$ satisf\/ies the condition in Def\/inition~\ref{1.6}.

Therefore take $\varepsilon>0$ and vectors $(\eta_k)$ in $\mathcal H$.
First choose an index $i_0$ so that
\begin{gather*}
\| \xi \otimes \eta_k - \xi_{i_0} \otimes \eta_k \| < \varepsilon
\qquad
\text{and}
\qquad
\| \xi' \otimes \eta_k - K\xi_{i_0} \otimes \eta_k \| < \varepsilon
\end{gather*}
for all~$k$.
Then choose the elements $(p_j)$ and $(q_j)$ as in Def\/inition~\ref{1.6} for the pair $(\xi_{i_0},K\xi_{i_0})$.
These elements will now also satisfy the inequalities for the original pair $(\xi,\xi')$, with $2\varepsilon$ instead
of~$\varepsilon$.
This is what we had to show.
\end{proof}

The counterpart of the other result for $S_0$, namely that $S_0(xy)=S_0(y)S_0(x)$ when $x,y\in \mathcal D_0$, is the
following.

\begin{Proposition}\label{1.10}
Let $x\in\mathcal D_0$ and assume that $x_1$ is as in Definition~{\rm \ref{1.2}}.
If $\xi\in\mathcal D(K)$ then $x\xi\in\mathcal D(K)$ and $Kx\xi=x_1K\xi$.
\end{Proposition}

\begin{proof}
Take a~pair $(x,x_1)$ of elements in~$M$ satisfying the condition as in Def\/inition~\ref{1.2}.
Take $\xi\in D(K)$ and put $\xi_1=K\xi$.
We must show that the pair $(x\xi,x_1\xi_1)$ satisf\/ies the conditions as in Def\/inition~\ref{1.6}.

To show this, take $\varepsilon>0$ and vectors $(\eta_k)$ in $\mathcal H$.
First choose elements $(p_i)$ and $(q_i)$ in~$M$ so that
\begin{gather*}
\Big \|x \xi\otimes \eta_k - \textstyle\sum\limits_i \Delta(p_i)(\xi\otimes q_i^*\eta_k)\Big\| < \varepsilon,
\qquad
 \Big\|x_1 \xi_1\otimes \eta_k - \textstyle\sum\limits_i \Delta(q_i)(\xi_1\otimes p_i^*\eta_k)\Big\| < \varepsilon
\end{gather*}
for all~$k$ as in Def\/inition~\ref{1.2}.
Next take $\varepsilon'>0$ and choose elements $(r_{ij})$ and $(s_{ij})$ in $\mathcal N_\psi$ so that
\begin{gather*}
 \Big\|\xi\otimes q_i^*\eta_k - V\Big(\textstyle\sum\limits_j \Lambda_\psi(r_{ij}) \otimes s_{ij}^*q_i^*\eta_k\Big)\Big\| < \varepsilon',
\\
 \Big\|\xi_1\otimes p_i^*\eta_k - V\Big(\textstyle\sum\limits_j \Lambda_\psi(s_{ij})\otimes r_{ij}^*p_i^*\eta_k\Big)\Big\| < \varepsilon'
\end{gather*}
for all~$i$ and all~$k$ as in Def\/inition~\ref{1.6}.
Then we f\/ind for all~$k$ on the one hand
\begin{gather*}
\Big\|x\xi\otimes\eta_k -  V\Big(\textstyle\sum\limits_{ij} \Lambda_\psi(p_ir_{ij}) \otimes s_{ij}^*q_i^*\eta_k\Big)\Big\|
\leq \|x \xi\otimes \eta_k - \textstyle\sum\limits_i \Delta(p_i)(\xi\otimes q_i^*\eta_k)\|
\\
\qquad
\phantom{\leq}{}
+ \Big\|\textstyle\sum\limits_i \Delta(p_i)(\xi\otimes q_i^*\eta_k)- V\Big(\textstyle\sum\limits_{ij} \Lambda_\psi(p_ir_{ij})
\otimes s_{ij}^*q_i^*\eta_k\Big)\Big\|
\\
\qquad{}
 \leq \varepsilon + \Big\|\textstyle\sum\limits_i \Delta(p_i)\Big((\xi\otimes q_i^*\eta_k)- V\Big(\textstyle\sum\limits_j
\Lambda_\psi(r_{ij}) \otimes s_{ij}^*q_i^*\eta_k\Big)\Big)\Big\|
\\
\qquad{}
 \leq \varepsilon + \textstyle\sum\limits_i\| p_i\| \Big\| \xi\otimes q_i^*\eta_k- V\Big(\textstyle\sum\limits_j
\Lambda_\psi(r_{ij}) \otimes s_{ij}^*q_i^*\eta_k\Big)\Big\|
 \leq \varepsilon + \varepsilon'\textstyle\sum\limits_i \|p_i\|.
\end{gather*}
Similarly on the other hand
\begin{gather*}
\Big\|x_1\xi_1\otimes\eta_k - V\Big(\textstyle\sum\limits_{ij} \Lambda_\psi(q_is_{ij}) \otimes r_{ij}^*p_i^*\eta_k\Big)\Big\|\leq
\varepsilon + \varepsilon' \textstyle\sum\limits_i \| q_i\|
\end{gather*}
for all~$k$.

If we choose $\varepsilon'$ so that $\varepsilon'\sum\limits_i \|p_i\|<\varepsilon$ and $\varepsilon'
\textstyle\sum\limits_i \| q_i\|<\varepsilon$, we can complete the proof.
\end{proof}

A possible proof of the formula $S_0(xy)^*=S_0(x)^*S_0(y)^*$ when $x,y\in\mathcal D_0$ would be of the same type as the one above.

As an {\it important consequence} of the above proposition we f\/ind that, if $\mathcal D(K)$ is dense and if~$x$ and
$x_1$ are as in Def\/inition~\ref{1.2}, then $x=0$ will imply $x_1=0$.
Indeed, it will follow that $x_1K\xi=0$ for all $\xi\in\mathcal K$ and by Proposition~\ref{1.9}i) we have that the
range of~$K$ is equal to $\mathcal D(K)$.

\medskip

{\bf Density of the domain of the operator~$\boldsymbol{K}$.}
For the {\it following step}, we need a~left Haar weight on~$M$.
It is used to produce (enough) elements in $\mathcal D(K)$ and in~$\mathcal D_0$.
Recall that a~{\it left Haar weight} is a~faithful normal semi-f\/inite weight~$\varphi$ on~$M$ satisfying left
invariance, i.e.\
\begin{gather*}
\varphi((\omega\otimes\iota)\Delta(x))=\omega(1)\varphi(x),
\end{gather*}
whenever $x\in M$, $x\geq 0$ and $\varphi(x)<\infty$ and when $\omega\in M_*$ and $\omega\geq 0$.
For a~left Haar weight we have the {\it left regular representation}.
We use $\mathcal H_\varphi$ for the G.N.S.-space and $\Lambda_\varphi:\mathcal N_\varphi \to \mathcal H_\varphi$ for the
associated canonical map.
Again we let~$M$ act directly on $\mathcal H_\varphi$ (i.e.\ we drop the notation $\pi_\varphi$ as we did before with~$\psi$).
Later we will identify the two Hilbert spaces $\mathcal H_\varphi$ and $\mathcal H_\psi$ in such a~way that the actions
of~$M$ are the same (see the end of Section~\ref{Section3}).

The left regular representation is considered in the next proposition.

\begin{Proposition}\label{1.11}
There is a~bounded operator~$W$ on $\mathcal H\otimes \mathcal H_\varphi$, characterized $($and defined$)$~by
\begin{gather*}
((\omega\otimes\iota)W^*)\Lambda_\varphi(x)=\Lambda_\varphi((\omega\otimes\iota)\Delta(x)),
\end{gather*}
when $x\in \mathcal N_\varphi$ and $\omega\in \mathcal B(\mathcal H)_*$.
Now, the `first leg' of~$W$ sits in~$M$, that is $W\in M\otimes \mathcal B(\mathcal H_\varphi)$ and we have:
\begin{enumerate}\itemsep=0pt
\item[$i)$] $WW^*=1$ $($i.e.~$W$ is a~co-isometry$)$,

\item[$ii)$] $(1\otimes x)W=W\Delta(x)$ for all $x\in M$,

\item[$iii)$] $(\Delta\otimes\iota)W=W_{13}W_{23}$.
\end{enumerate}
\end{Proposition}

Here, roughly speaking, we have $W^*(\xi\otimes\Lambda_\varphi(x)){=}\sum\limits_{(x)} x_{(1)}\xi \otimes
\Lambda_\varphi(x_{(2)})$ when $\Delta(x){=}\sum\limits_{(x)} x_{(1)}\otimes x_{(2)}$ formally.
Observe the dif\/ference in convention (using the adjoint) when compared with the right regular representation (cf.\ Proposition~\ref{1.5}).
The proof of this proposition however is completely similar as for the right regular representation.

In order to use~$W$ to construct elements in $\mathcal D(K)$ and in $\mathcal D_0$, we need dif\/ferent steps.
We formulate dif\/ferent lemmas as some of the results will be needed later.
First we have the following.

\begin{Lemma}\label{1.12}
Let $\omega \in \mathcal B(\mathcal H_\varphi)_*$ and $x=(\iota\otimes\omega)W$ and
$x_1=(\iota\otimes\overline\omega)W$, then $x\in \mathcal D_0$ and $x_1$ satisfies the conditions as in
Definition~{\rm \ref{1.2}}.
\end{Lemma}

\begin{proof}
Assume that $\omega=\langle\,\cdot\, \xi,\eta\rangle$.
Take an orthonormal basis $(\xi_j)$ in $\mathcal H_\varphi$.
Def\/ine
\begin{gather*}
p_j=(\iota\otimes \langle\,\cdot\, \xi_j,\eta\rangle)W
\qquad
\text{and}
\qquad
q_j=(\iota\otimes \langle\,\cdot\, \xi_j,\xi\rangle)W.
\end{gather*}
Using the formula $(\Delta\otimes\iota)W=W_{13}W_{23}$ (Proposition~\ref{1.11}), we f\/ind
\begin{gather*}
\textstyle\sum\limits_j \Delta(p_j)(1\otimes q_j^*) = (\iota\otimes\iota\otimes\omega)(((\Delta\otimes\iota)W)(1\otimes W^*))
\\
\phantom{\textstyle\sum\limits_j \Delta(p_j)(1\otimes q_j^*)}
 = (\iota\otimes\iota\otimes\omega)(W_{13}W_{23}W_{23}^*)
 = (\iota\otimes\iota\otimes\omega)W_{13}
 = x\otimes 1
\end{gather*}
and similarly
\begin{gather*}
\textstyle\sum\limits_j \Delta(q_j)(1\otimes p_j^*) = (\iota\otimes\iota\otimes\overline\omega)(W_{13}W_{23}W_{23}^*)
 = (\iota\otimes\iota\otimes\overline\omega)W_{13}
 = x_1\otimes 1.
\end{gather*}
The sums converge in the strong operator topology.

This gives the result for elements~$\omega$ of the form $\langle\,\cdot\, \xi,\eta\rangle$.
Then it follows for all $\omega\in \mathcal B(\mathcal H_\varphi)_*$ by approximation.
\end{proof}

Remark that only the essential properties of~$W$ are used in the above argument and that it is not necessary to have
a~left regular representation, associated to a~left Haar weight.
Only the conditions i) and iii) of Proposition~\ref{1.11} are needed.

Compare this lemma with Proposition~5.6 in~\cite{V-VD} where a~similar argument is found.
Observe again that one of the dif\/ferences between this approach to the antipode and the one in~\cite{V-VD} lies in the
fact that we avoid the use of operator space techniques here.

Later we will combine this result with the property proven in Proposition~\ref{1.10} (cf.\ Proposition~\ref{1.16}).

In a~similar way, elements in the domain of~$K$ are constructed, but here we have to be a~bit more careful.
First we have the following lemma.

\begin{Lemma}\label{1.13}
If $c\in \mathcal N_\psi$ and $\omega\in \mathcal B(\mathcal H_\varphi)_*$ we have $(\iota\otimes\omega(c\,\cdot\,))W\in\mathcal N_\psi$.
\end{Lemma}

\begin{proof}
If we let $x= (\iota\otimes\omega(c\,\cdot\,))W$, we get
\begin{gather*}
x^*x \leq \|\omega\| (\iota\otimes|\omega|)(W^*(1\otimes c^*)(1\otimes c)W)
 = \|\omega\| (\iota\otimes|\omega|)(\Delta(c^*)W^*W\Delta(c))
\\
\phantom{x^*x}
 \leq \|\omega\| (\iota\otimes|\omega|)(\Delta(c^*)\Delta(c))
 = \|\omega\| (\iota\otimes|\omega|)(\Delta(c^*c)).
\end{gather*}
As~$\psi$ is right invariant and $c\in \mathcal N_\psi$, we get also $x\in \mathcal N_\psi$.
\end{proof}

Observe that we do not need that~$W$ is unitary.
It is suf\/f\/icient for this argument that $W^*W\leq 1$ and this is true for a~co-isometry.

Now the following result should not come as a~surprise.

\begin{Lemma}\label{1.14}
Let $c,d\in \mathcal N_\psi$ and $\omega\in \mathcal B(\mathcal H_\varphi)_*$ and define
$\xi=\Lambda_\psi((\iota\otimes\omega(c \cdot d^*))W)$.
Then $\xi\in \mathcal D(K)$ and $K\xi=\Lambda_\psi((\iota\otimes\overline\omega(d \cdot c^*))W)$.
\end{Lemma}

\begin{proof}
The proof of this lemma is based on the same decomposition as in Lemma~\ref{1.12}.

Take~$\omega$ of the form $\langle\,\cdot\,\xi',\eta'\rangle$ where $\xi'$ and $\eta'$ are vectors in $\mathcal H_\varphi$.
Take an orthonormal basis~$(\xi_j)$ in~$\mathcal H_\varphi$.
Def\/ine elements in~$M$ as before~by{\samepage
\begin{gather*}
p_j =(\iota\otimes\langle \,\cdot\,\xi_j,c^*\eta'\rangle)W,
\qquad
q_j =(\iota\otimes\langle \,\cdot\,\xi_j,d^*\xi'\rangle)W.
\end{gather*}
By Lemma~\ref{1.13} we have $p_j, q_j\in \mathcal N_\psi$.
This is necessary for the use of Def\/inition~\ref{1.6}.}

We know that
\begin{gather*}
\textstyle\sum\limits_j \Delta(p_j)(1\otimes q_j^*)  = x\otimes 1,
\qquad
\textstyle\sum\limits_j \Delta(q_j)(1\otimes p_j^*)  = x_1\otimes 1,
\end{gather*}
where
\begin{gather*}
x =(\iota\otimes\langle\,\cdot\, d^*\xi',c^*\eta'\rangle)W,
\qquad
x_1 =(\iota\otimes\langle\,\cdot\, c^*\eta',d^*\xi'\rangle)W
\end{gather*}
as in Lemma~\ref{1.12}, with convergence in in the strong operator topology.
Because now all these elements belong $\mathcal N_\psi$, using the properties of the map $\Lambda_\psi$, we will also
have
\begin{gather*}
V\Big(\textstyle\sum\limits_j \Lambda_\psi(p_j)\otimes q_j^*\eta\Big)  = \Lambda_\psi(x)\otimes \eta,
\qquad
V\Big(\textstyle\sum\limits_j \Lambda_\psi(q_j)\otimes p_j^*\eta\Big)  = \Lambda_\psi(x_1)\otimes \eta
\end{gather*}
for all $\eta\in \mathcal H$.
Now convergence will be in the norm topology of the Hilbert space tensor product.
Then it follows from Def\/inition~\ref{1.6} that $\Lambda_\psi(x)\in \mathcal D(K)$ and that
$K\Lambda_\psi(x)=\Lambda_\psi(x_1)$.
This is what we had to show.
\end{proof}

Having these results, we are ready to show that the domain of~$K$ is dense.
Simultaneously, we obtain that the right regular representation~$V$ is unitary.
Indeed, as the proof of the two results are intimately related, we formulate them below in one proposition.

\begin{Proposition}\label{1.15}
The operator~$V$ is unitary.
And the operator~$K$ is densely defined.
\end{Proposition}

\begin{proof}
For the proof of the f\/irst statement, we use Kustermans' trick as in~\cite{K-V2}.
Def\/ine
\begin{gather*}
\mathcal K= \overline{\SP} \big\{\Lambda_\psi((\iota\otimes\omega(c\,\cdot\,))W) \,|\, c\in \mathcal N_\psi,\; \omega\in
\mathcal B(\mathcal H_\varphi)_*\big\},
\end{gather*}
where by $\overline{\SP}$ we mean the closed linear span.

Consider~$V$ as acting on the space $\mathcal H_\psi\otimes \mathcal H_\varphi$ by taking $\mathcal H_\varphi$ for $\mathcal H$.

Consider the notations of Lemma~\ref{1.14}, but with $d=1$.
In this case
$\xi=\Lambda_\psi((\iota\otimes\omega(c\,\cdot\,))W)$.
Using the same techniques as in the proof of Lemmas~\ref{1.12} and~\ref{1.14}, we f\/ind that $\xi\otimes \eta$ is
approximated by f\/inite sums of the form $\sum V(\Lambda_\psi(p_j)\otimes q_j^*\eta)$ for any $\eta\in \mathcal H_\varphi$.
Because $\xi$, as well as all the elements $\Lambda_\psi(p_j)$ belong to $\mathcal K$,
we f\/ind that $\mathcal K\otimes\mathcal H_\varphi \subseteq V(\mathcal K\otimes \mathcal H_\varphi)$.

On the other hand, we have the formula
\begin{gather*}
(\iota\otimes \varphi)((\Delta(x^*)(1\otimes y))=(\iota\otimes\langle\,\cdot\,\Lambda_\varphi(y),\Lambda_\varphi(x)\rangle)W,
\end{gather*}
whenever $x,y\in \mathcal N_\varphi$.
This follows easily from the def\/ining formula for~$W$ in Proposition~\ref{1.11}.
We can now approximate any linear functional of the form $\omega(c\,\cdot\,)$ by functionals of the form
$\langle\,\cdot\,\Lambda_\varphi(y),\Lambda_\varphi(x)\rangle$, where we take $y\in \mathcal N_\varphi$ and $x\in
\mathcal N_\varphi\cap \mathcal N_\psi^*$.
Moreover, we can approximate any element in $M_*$ by linear functionals of the form $\varphi(\,\cdot\,y)$ with
appropriate elements $y\in \mathcal N_\varphi$.
As a~consequence of all these carefully chosen approximations, we f\/ind that also
\begin{gather*}
\mathcal K= \overline{\SP} \big\{\Lambda_\psi((\iota\otimes\omega)\Delta(x)) \,|\, x\in \mathcal N_\psi,\; \omega\in
M_* \big\}.
\end{gather*}
Consequently we see that also $V(\mathcal H_\psi\otimes\mathcal H_\varphi)\subseteq \mathcal K\otimes \mathcal H_\varphi$.

By a~combination of the two results above and using that~$V$ is isometric, we get $\mathcal K=\mathcal H_\psi$.
Therefore~$V$ is unitary.

As we also have
\begin{gather*}
\mathcal K= \overline{\SP} \big\{\Lambda_\psi((\iota\otimes\omega(c\,\cdot\,d^*))W) \,|\, c,d\in \mathcal N_\psi,\;
\omega\in \mathcal B(\mathcal H_\varphi)_* \big\},
\end{gather*}
we get from $\mathcal K=\mathcal H_\psi$ that $\mathcal D(K)$ is dense.
\end{proof}

Compare the proof of this proposition with arguments found in~\cite[Section~3.3]{K-V2}.

By symmetry, of course also the left regular representation~$W$ associated to any left Haar weight will be unitary.
Observe that we now can rewrite the formulas~ii) of Proposition~\ref{1.5} and~ii) of Proposition~\ref{1.11} as
$\Delta(x)=V(x\otimes 1)V^*$ and $\Delta(x)=W^*(1\otimes x)W$ respectively.

The unitarity of the regular representations can also be proven in an other, perhaps shorter (still essentially the
same) way, but because we also need the density of the domain $\mathcal D(K)$ of~$K$, we have chosen to prove these
results together as above.

It should not come as a~surprise that the density of $\mathcal D(K)$ is essentially the same result as saying that the
isometry~$V$ is in fact a~unitary.
Indeed, when $\mathcal D(K)$ is dense, it follows that for any $\xi\in\mathcal H_\psi$ and any $\eta\in\mathcal H$, the
vector $\xi\otimes\eta$ can be approximated with elements in the range of~$V$ (cf.\ Def\/inition~\ref{1.6}).
Roughly speaking, this says that the map $p\otimes q \mapsto \Delta(p)(1\otimes q)$, considered on the Hilbert space
level, has dense range.
Later, at the end of this section, we will see that this map on $M\otimes M$ also has dense range.
This in turn will be a~consequence of the density of $\mathcal D_0$ (cf.\ Proposition~\ref{1.21} below).
Remark that, although there are similarities, the two density results are dif\/ferent because the topologies considered are dif\/ferent.

\medskip

{\bf The antipode and its polar decomposition.}
We have now shown that the domain $\mathcal D(K)$ of the operator~$K$
is dense and as we mentioned already (see the remark following Proposition~\ref{1.10}) it would now be possible to
def\/ine the antipode as the map $S_0$ given by $S_0(x)=x_1^*$ (cf.\ the~remark following Def\/inition~\ref{1.2}).
It would still be necessary to show that also the domain~$\mathcal D_0$ of~$S_0$ is dense.

Eventually we will see that indeed, the space $\mathcal D_0$ is dense (see Proposition~\ref{1.21} below).
But f\/irst we will construct the antipode by means of its polar decomposition.
Some of the formulas needed to do this will also play an important role in the next section where we obtain the main results.

Let us f\/irst formulate a~result that easily follows from combining Lemma~\ref{1.12} with Proposition~\ref{1.10}.

\begin{Proposition}\label{1.16}
For any $\xi\in \mathcal D(K)$ and $\omega\in\mathcal B(\mathcal H_\varphi)_*$, we have that
$((\iota\otimes\omega)W)\xi\in \mathcal D(K)$ and
\begin{gather*}
K((\iota\otimes\omega)W)\xi=((\iota\otimes\overline\omega)W)K\xi.
\end{gather*}
\end{Proposition}

Next we need a~similar formula, but for the other leg of~$W$.
And because eventually we will need all of this to obtain uniqueness of the Haar weights, we will work with two left
Haar weights~$\varphi_1$ and~$\varphi_2$.
We will use the left regular representations for these two left Haar weights and we will use $W_1$ and $W_2$ to denote them.
We will in what follows consider these operators as acting on the spaces $\mathcal H_\psi\otimes\mathcal H_{\varphi_1}$
and $\mathcal H_\psi\otimes\mathcal H_{\varphi_2}$ respectively.

We then have the following result.

\begin{Proposition}\label{1.17}
Let $T_r$ be the closure of the operator $\Lambda_{\varphi_1}(x)\mapsto \Lambda_{\varphi_2}(x^*)$ with $x\in \mathcal
N_{\varphi_1} \cap \mathcal N_{\varphi_2}^*$.
If $\xi\in \mathcal D(T_r)$ then $((\omega\otimes\iota)W_1^*)\xi \in\mathcal D(T_r)$ for all $\omega\in\mathcal
B(\mathcal H_\psi)_*$ and
\begin{gather*}
T_r((\omega\otimes\iota)W_1^*)\xi= ((\overline\omega\otimes\iota)W_2^*)T_r\xi.
\end{gather*}
\end{Proposition}

\begin{proof}
Fix $\omega\in \mathcal B(\mathcal H_\psi)_*$.
First we prove the formula for $\xi=\Lambda_{\varphi_1}(x)$ with $x\in \mathcal N_{\varphi_1} \cap \mathcal
N_{\varphi_2}^*$.
We get
\begin{gather*}
((\omega\otimes\iota)W_1^*)\Lambda_{\varphi_1}(x) = \Lambda_{\varphi_1}((\omega\otimes\iota)\Delta(x))
\end{gather*}
by the def\/inition of $W_1$, see Proposition~\ref{1.11}.
Then, because also $(\omega\otimes\iota)\Delta(x)\in \mathcal N_{\varphi_1} \cap \mathcal N_{\varphi_2}^*$, we get from
the def\/inition of $T_r$ that
\begin{gather*}
\begin{split}
& T_r(((\omega\otimes\iota)W_1^*)\Lambda_{\varphi_1}(x)) = \Lambda_{\varphi_2}((\overline\omega\otimes\iota)\Delta(x^*)) = ((\overline\omega\otimes\iota)W_2^*)\Lambda_{\varphi_2}(x^*)\\
& \hphantom{T_r(((\omega\otimes\iota)W_1^*)\Lambda_{\varphi_1}(x))}{}
 = ((\overline\omega\otimes\iota)W_2^*)T_r\Lambda_{\varphi_1}(x).
\end{split}
\end{gather*}
The result for any vector $\xi\in \mathcal D(T_r)$ follows because $T_r$ is the closure of the map
$\Lambda_{\varphi_1}(x)\mapsto \Lambda_{\varphi_2}(x^*)$ with $x\in \mathcal N_{\varphi_1} \cap \mathcal N_{\varphi_2}^*$.
\end{proof}

Now we will combine the above result with the similar formula for~$K$, applied for both $W_1$ and $W_2$.
Compare with results in~\cite[Section~5.2]{K-V2}.

\begin{Proposition}\label{1.18}
With the notations as before, we have the equality
\begin{gather*}
(K\otimes T_r)W_1 = W_2^*(K\otimes T_r).
\end{gather*}
\end{Proposition}

\begin{proof}
Take vectors $\xi\in \mathcal D(K)$, $\xi'\in \mathcal D(K^*)$, $\eta\in \mathcal D(T_r)$ and $\eta'\in \mathcal D(T_r^*)$.
Remember that $\xi, \xi'\in \mathcal H_\psi$ while $\eta\in \mathcal H_{\varphi_1}$ and $\eta'\in \mathcal H_{\varphi_2}$.

i) We f\/irst use the formula with~$K$ (as proven in Proposition~\ref{1.16}).
Then we f\/ind
\begin{gather*}
\langle W_2(K\xi \otimes T_r\eta), \xi'\otimes \eta'\rangle
 =\langle ((\iota\otimes\langle\,\cdot\,T_r\eta,\eta'\rangle)W_2)K\xi, \xi'\rangle
 =\langle K((\iota\otimes\langle\,\cdot\,\eta',T_r\eta\rangle)W_2)\xi, \xi'\rangle
\\
\phantom{\langle W_2(K\xi \otimes T_r\eta), \xi'\otimes \eta'\rangle}
 =\langle ((\iota\otimes\langle\,\cdot\,\eta',T_r\eta\rangle)W_2)\xi, K^*\xi'\rangle^-
 =\langle W_2(\xi\otimes\eta'),K^*\xi'\otimes T_r\eta\rangle^-.
\end{gather*}
Next we write this last expression as
\begin{gather*}
\langle((\langle\,\cdot\, \xi, K^*\xi'\rangle\otimes \iota)W_2)\eta', T_r\eta\rangle^-.
\end{gather*}
If we now use the formula with $T_r$ from Proposition~\ref{1.17}, we f\/ind by a~similar calculation that this is equal to
\begin{gather*}
\langle W_1^*(\xi\otimes\eta),K^*\xi'\otimes T_r^*\eta'\rangle^-.
\end{gather*}
This implies the inclusion $ W_2(K\otimes T_r) \subseteq (K\otimes T_r)W_1^*$.

ii) On the other hand, if we proceed as
above, but now f\/irst using the formula for $T_r$ and then the formula for~$K$, we f\/ind
\begin{gather*}
\langle W_2^*(K\xi \otimes T_r\eta), \xi'\otimes \eta'\rangle = \langle W_1(\xi\otimes\eta),K^*\xi'\otimes
T_r^*\eta'\rangle^-.
\end{gather*}
This in turn implies the inclusion $ W_2^*(K\otimes T_r) \subseteq (K\otimes T_r)W_1$.

If we take the f\/irst inclusion and apply $W_2^*$ from the left and $W_1$ from the right, we get $ (K\otimes T_r)W_1
\subseteq W_2^*(K\otimes T_r)$.
If we combine this with the previous inclusion, we get the result.
\end{proof}

This formula is very important for the further development in Section~\ref{Section3}.
Of course, we can also replace both $W_1$ and $W_2$ by~$W$ associated with any left Haar weight~$\varphi$.
We will use both cases in the next section when we show that Haar weights are unique.
In this section, we will use it (with~$W$) to construct the antipode and to prove some more density results as we announced.

In order to use our formula, we need to consider the polar decomposition of the operators involved.

\begin{Notation}\label{1.19}
Let~$K$ be the operator on $\mathcal H_\psi$ as def\/ined in Def\/initions~\ref{1.6} and~\ref{1.8}. Now let~$T$ be
the closure of the map $\Lambda_\varphi(x)\mapsto \Lambda_\varphi(x^*)$ where $x\in \mathcal N_\varphi \cap \mathcal
N_\varphi^*$.
We use
\begin{gather*}
K=IL^{\frac12}
\qquad
\text{and}
\qquad
T=J\nabla^{\frac12}
\end{gather*}
to denote the polar decompositions of these operators.
\end{Notation}

The properties of all these operators are well-known and easy consequences of the fact that~$K$ and~$T$ are conjugate
linear and involutive.
We have e.g.\ that $J\nabla J =\nabla^{-1}$ so that $J\nabla^{it} J= \nabla^{it}$ (because~$J$ is conjugate linear).
Similarly for the operators~$I$ and~$L$.
See e.g.\ Chapter VI in~\cite{T3}.

Remember that, roughly speaking, our operator~$K$ coincides with the operator $G^*$ in~\cite{K-V2} and therefore, that
the operator~$L$ is essentially the operator $N^{-1}$ in~\cite[Section~5]{K-V2}.
See also Section~\ref{Section5}, in particular Remark~\ref{4.10}.

If we apply Proposition~\ref{1.18} to the case $\varphi_1=\varphi_2=\varphi$, we get $(K\otimes T)W=W^*(K\otimes T)$
where~$W$ is the left regular representation associated with~$\varphi$.
As a~consequence of the uniqueness of the polar decomposition, we get the following result.

\begin{Proposition}\label{1.20}
We have $(I\otimes J)W(I\otimes J)=W^*$ and also $(L^{it}\otimes \nabla^{it})W(L^{-it}\otimes \nabla^{-it})=W$ for all
$t\in {\mathbb R}$.
\end{Proposition}

In the next section, we will use similar formulas, but for two weights and we will combine them with these formulas here
to get uniqueness of the Haar weights.

We will now show in the next proposition that the left leg of~$W$ is dense in~$M$.
Therefore, the above formulas will allow us to def\/ine maps $R:M\to M$ and $\tau_t:M\to M$ for all~$t$ by $R(x)=I x^* I$
and $\tau_t(x)=L^{it}xL^{-it}$.
These maps will give us the polar decomposition of the antipode (see Def\/inition~\ref{1.23} below).

We f\/irst need the following observation.
Denote by $(\sigma_t^\varphi)_{t\in {\mathbb R}}$ the modular automorphisms on~$M$ def\/ined by
$\sigma_t^\varphi(x)=\nabla^{it}x\nabla^{-it}$.
Similarly let us def\/ine the one-parameter group of automor\-phisms~$(\tau_t)$ on $\mathcal B(\mathcal H_\psi)$~by
$\tau_t=L^{it}\,\cdot\,L^{-it}$.
Then it follows from the second formula in Proposition~\ref{1.20} and from $\Delta(x)=W^*(1\otimes x)W$ for all $x\in M$
that $\Delta(\sigma_t^\varphi(x)))=(\tau_t\otimes\sigma_t^\varphi)\Delta(x)$ for all $x\in M$.
From this formula, it follows that the space of slices, spanned by the elements $(\omega\otimes\iota)\Delta(x)$ with
$x\in M$ and $\omega\in M_*$ will be left invariant under the modular automorphisms $(\sigma_t^\varphi)_{t\in {\mathbb R}}$.
We will need this for the proof of the following proposition (see in~\cite[Proposition~1.4]{K-V2}).

\begin{Proposition}\label{1.21}
Let~$W$ be the left regular representation associated with some left Haar weight~$\varphi$ as before.
Then the following three subspaces of~$M$
\begin{enumerate}\itemsep=0pt
\item[$i)$] $\{(\iota\otimes\omega)W \,|\, \omega\in \mathcal
B(\mathcal H_\varphi)_*\}$,

\item[$ii)$] $\SP \{(\iota\otimes\omega)\Delta(x) \,|\, x\in M,\; \omega\in M_* \}$,

\item[$iii)$] $\SP \{(\omega\otimes\iota)\Delta(x) \,|\, x\in M,\; \omega\in M_* \}$
\end{enumerate}
are~$\sigma$-weakly dense in~$M$.
\end{Proposition}

\begin{proof}
We will only consider i) and ii) because the density in iii) will follow by symmetry.

We f\/irst claim that the spaces in i) and in ii) have the same closure.
This follows from the formula
\begin{gather*}
(\iota\otimes \varphi)((\Delta(x^*)(1\otimes y))=(\iota\otimes\langle\,\cdot\,\Lambda_\varphi(y),\Lambda_\varphi(x)\rangle)W
\end{gather*}
with $x,y\in \mathcal N_\varphi$ considered already in the proof of Proposition~\ref{1.15}.

Let us now denote by $M_e$ the closure of the space in i).
We have to show that this is equal to~$M$.

It follows from the fact that~$W$ satisf\/ies the pentagon equation, that~$M_e$ is a~subalgebra of~$M$.
Because~$M_e$ is also the closure of the space in ii) and this is obviously self-adjoint, we get that~$M_e$ is
a~$*$-subalgebra of~$M$.

In the proof of Proposition~\ref{1.15}, we have seen that the space
\begin{gather*}
\SP \big\{\Lambda_\psi((\iota\otimes\omega)\Delta(x)) \,|\, x\in\mathcal N_\psi,\; \omega\in M_* \big\}
\end{gather*}
is dense in $\mathcal H_\psi$.
Standard approximation techniques give that also
\begin{gather*}
\SP \big\{\Lambda_\psi((\iota\otimes\omega)\Delta(x)) \,|\, x\in \mathcal N_\psi\cap \mathcal N_\psi^*,\; \omega\in M_*\big\}
\end{gather*}
is still dense in $\mathcal H_\psi$.
This will imply that the space $\Lambda_\psi(\mathcal N_\psi\cap \mathcal N_\psi^*\cap M_e)$ is dense in
$\Lambda_\psi(\mathcal N_\psi\cap \mathcal N_\psi^*)$.

We have seen in a~remark preceding this proposition, that the space of slices in iii) is inva\-riant by the modular automorphisms of~$\varphi$.
Similarly, the modular automorphism $\sigma_t^\psi$ leaves~$M_e$ invariant.
So the space $\Lambda_\psi(\mathcal N_\psi\cap \mathcal N_\psi^*\cap M_e)$ will be invariant under the modular unitaries
$\nablad^{it}$ (where we use $\nablad$ for the modular operator associated with the right Haar weight~$\psi$).
It follows that also $\Lambda_\psi(\mathcal N_\psi\cap \mathcal N_\psi^*\cap M_e)$ is dense in $\Lambda_\psi(\mathcal
N_\psi\cap \mathcal N_\psi^*)$ with respect to the $^\#$-norm (cf.~Section~1 in Chapter~VI in~\cite{T3}).
Then, from a~result in Hilbert algebra theory (see Lemma~5.1 and the proof of Theorem~10.1 in~\cite{T1}),
it will follow that also $\mathcal N_\psi\cap \mathcal N_\psi^*\cap M_e$ is dense in~$M$.
Therefore we have that $M_e=M$.
This completes the proof of the proposition.
\end{proof}

From this proposition and taking into account Lemma~\ref{1.12}, it follows that the map $p\otimes q\mapsto
\Delta(p)(1\otimes q)$ has dense range in $M\otimes M$.
Indeed, one can approximate elements of the form $x\otimes 1$ by linear combinations of elements of the form
$\Delta(p)(1\otimes q)$ when $x\in \mathcal D_0$ (by the very def\/inition of $\mathcal D_0$) and as now this domain is
shown to be dense, we can do this for any $x\in M$.
By multiplying with elements of~$M$ in the second factor, we get the density of $\Delta(M)(1\otimes M)$.
By symmetry, we also have that $\Delta(M)(M\otimes 1)$ will be dense in $M\otimes M$.
Also compare with the remark following Proposition~\ref{1.15}.

We {\it finish this section} by the def\/inition of the antipode~$S$ and its polar decomposition and by formulating an
important property which will be frequently used in the next sections.

\begin{Definition}\label{1.22}
Def\/ine $R:M\to M$ by $R(x)=Ix^*I$ and $\tau_t:M\to M$ by $\tau_t(x)=L^{it}xL^{-it}$.
\end{Definition}

It is a~consequence of the density result i) in Proposition~\ref{1.20} and Proposition~\ref{1.21} that these maps do
leave~$M$ invariant.
We have that~$R$ and $\tau_t$ commute because~$I$ and $L^{it}$ commute for all~$t$.
The automorphisms $(\tau_t)$ are called {\it the scaling automorphisms} whereas the anti-automorphism~$R$ is called the
{\it unitary antipode}.
Together they represent what is commonly referred to as {\it the polar decomposition of the antipode}:

\begin{Definition}\label{1.23}
The antipode~$S$ is def\/ined as the composition $R \tau_{-\frac{i}{2}}$ where $\tau_{-\frac{i}{2}}$ is the analytic
generator associated with the one-parameter group $(\tau_t)$ in the point $-\frac{i}{2}$.
\end{Definition}

Recall, as was already mentioned in the introduction, that the analytic generator is def\/ined f\/irst on the predual $M_*$
and then on~$M$ by taking the adjoint.

In the next sections, we will (essentially) no longer need the operators~$K$,~$I$ and~$L$ (they will be replaced~by
other, more adapted operators, see Section~\ref{Section5}).
We will use the unitary antipode~$R$ and the scaling group $(\tau_t)$.
The main results in this section involving~$K$, $I$ and~$L$ can be restated solely in terms of~$R$ and $\tau_t$.
This is e.g.\ quite obvious for the two formulas in Proposition~\ref{1.20} (and for the similar result involving two
left Haar weights, see a~remark in the beginning of Section~\ref{Section3}, before Proposition~\ref{2.2}).

It is somewhat harder with Proposition~\ref{1.16}, but the following result can be shown about the antipode as def\/ined in Def\/inition~\ref{1.23}.

\begin{Proposition}\label{1.24}
For any $\omega\in\mathcal B(\mathcal H_\varphi)_*$ we have that $(\iota\otimes \omega)W\in \mathcal D(S)$ and
$S((\iota\otimes\omega)W)=(\iota\otimes\omega)(W^*)$.
The space of such elements is invariant under the scaling automorphisms $(\tau_t)$ and it is a~core for~$S$.
\end{Proposition}

Indeed, formally, from Proposition~\ref{1.16}, we get $KxK=x_1^*$ when $x=(\iota\otimes\omega)W$ and
$x_1=(\iota\otimes\omega)(W^*)$.
Now
\begin{gather*}
KxK=I L^\frac12 x L^{-\frac12} I = R\big(\tau_{-\frac{i}{2}}(x)\big)^*.
\end{gather*}
One has to be somewhat careful, but the argument can be made precise.
Furthermore, because of the second formula in Proposition~\ref{1.20}, the space of such elements $(\iota\otimes\omega)W$
will be invariant under the scaling automorphisms and as this space is dense in~$M$, it will be a~core for~$S$.
See e.g.~\cite{VD8} for details.

We will refer to this result when in the sequel, we write loosely $(S\otimes \iota)W=W^*$.

\section[The main results about $(M,\Delta)$]{The main results about $\boldsymbol{(M,\Delta)}$}\label{Section3}

In Section~\ref{Section2}, we have introduced the antipode and obtained some results about density.
We also proved an important formula and some consequences of it.
In this section, we will prove the main results about a~locally compact quantum group and its related objects (such as
the left and right Haar weights, the modular automorphism groups, the scaling group and the unitary antipode, \dots).
In the next section, we will treat the dual and in the fourth section we will prove the main results about the objects
associated with a~locally compact quantum group in relation with those of the dual.

But f\/irst it is appropriate to formulate here the precise def\/inition of a~locally compact quantum group (in the von
Neumann algebraic setting), see~\cite{K-V3}.

\begin{Definition}\label{2.1}
Let~$M$ be a~von Neumann algebra.
Let~$\Delta$ be a~comultiplication on~$M$.
The pair $(M,\Delta)$ is called a~{\it locally compact quantum group} if there exist a~left and a~right Haar weight.
\end{Definition}

The correct def\/inition of a~comultiplication on a~von Neumann algebra~$M$ was given in Def\/inition~\ref{1.1}.
It is a~unital normal $*$-homomorphism from~$M$ to $M\otimes M$ satisfying coassociativity.
A~right Haar weight is a~faithful normal semi-f\/inite weight on~$M$ that is right invariant (see Def\/inition~\ref{1.4}).
A~left Haar weight is a~faithful normal semi-f\/inite weight on~$M$ that is left invariant.

In this section, $(M,\Delta)$ will be a~locally compact quantum group (in the sense of the above def\/inition) and~$\psi$
will be a~right Haar weight and~$\varphi$,~$\varphi_1$ and~$\varphi_2$ will denote left Haar weights.
We will use the notations and results of the previous section.

A~f\/irst major objective is to prove that Haar weights are unique.
We consider the two left Haar weights $\varphi_1$ and $\varphi_2$.
We denote by $(u_t)_{t\in {\mathbb R}}$ the Connes' cocycle Radon--Nikodym derivative of $\varphi_1$ w.r.t.\ $\varphi_2$ and
we will show that $u_t$ is a~scalar multiple of $1$ for each~$t$.
This will imply that the two weights are proportional.

Our starting point will be the formula in Proposition~\ref{1.18}.
Whereas in the previous section, we have used this formula for the case of one invariant weight, now we will also use it
for the two dif\/ferent weights $\varphi_1$ and $\varphi_2$.

Recall that we use $W_1$ and $W_2$ for the left regular representations associated with~$\varphi_1$ and~$\varphi_2$
respectively.
And as in the previous section, we denote by $T_r$ the closure of the operator $\Lambda_{\varphi_1}(x)\mapsto
\Lambda_{\varphi_2}(x^*)$ where $x\in \mathcal N_{\varphi_1} \cap \mathcal N_{\varphi_2}^*$.
Let $T_r=J_r\nabla_r^{\frac12}$ be the polar decomposition.
Remark that $\nabla_r$ is an operator on $\mathcal H_{\varphi_1}$ and that $J_r$ maps $\mathcal H_{\varphi_1}$ to
$\mathcal H_{\varphi_2}$.
Recall that we use $K=IL^{\frac12}$ for the polar decomposition of~$K$.
From Proposition~\ref{1.18}, we know that $(K\otimes T_r)W_1=W^*_2(K\otimes T_r)$.
Then, as in Proposition~\ref{1.20}, it follows from the uniqueness of the polar decompositions that $L^{it}\otimes
\nabla_r^{it}$ commutes with $W_1$ for all~$t$.
This result can be restated without the use of the operator~$L$.
We simply write e.g.\ $(\tau_t\otimes\iota)W_1=(1\otimes\nabla_r^{-it})W_1(1\otimes\nabla_r^{it})$.
We refer to the~remark following Def\/inition~\ref{1.23}.

Let us also use $T_1$ for the closure of the operator $\Lambda_{\varphi_1}(x)\mapsto \Lambda_{\varphi_1}(x^*)$ with $x
\in \mathcal N_{\varphi_1}\cap \mathcal N_{\varphi_1}^*$ and $T_1=J_1 \nabla_1^\frac12$ for its polar decomposition.
It is known from the {\it relative modular theory} (see e.g.\ Section~3 in Chapter~VIII of~\cite{T3}) that
$u_t=\nabla_1^{it}\nabla_r^{-it}$.
Then the following result follows easily.

\begin{Proposition}\label{2.2}
For each $t\in {\mathbb R}$ we have $\Delta(u_t)=1\otimes u_t$.
\end{Proposition}

\begin{proof}
We just saw that $(L^{it}\otimes \nabla_r^{it})W_1(L^{-it}\otimes \nabla_r^{-it})=W_1$ for all $t\in {\mathbb R}$.
If we combine this with the formula ii) in Proposition~\ref{1.20} of the previous section (for $W_1$), we get $(1\otimes
u_t)W_1(1\otimes u_t^*)=W_1$.
Because $W_1^*(1\otimes x)W_1=\Delta(x)$ for all $x\in M$, we get the result.
\end{proof}

If the right Haar weight~$\psi$ is bounded, we can immediately apply it on the formula above and get
$\psi(u_t)1=\psi(1)u_t$ so that $u_t$ must be a~scalar multiple of the identity for all~$t$.
We will be able to conclude this, also in general, but we need a~f\/iner argument.

Before we complete this argument, recall that we have def\/ined $R:M\to M$ by $R(x)=Ix^*I$ and $\tau_t:M\to M$ by
$\tau_t(x)=L^{it}xL^{-it}$ (see Def\/inition~\ref{1.22}).
Again we will also use $(\sigma_t^\varphi)_{t\in {\mathbb R}}$ to denote the modular automorphisms on~$M$ def\/ined~by
$\sigma_t^\varphi(x)=\nabla^{it}x\nabla^{-it}$ for the given left Haar weight~$\varphi$.

We get a~f\/irst set of important formulas.

\begin{Theorem}\label{2.3}
For all $x\in M$ and $t\in {\mathbb R}$ we have
\begin{enumerate}\itemsep=0pt
\item[$i)$] $\Delta(\sigma_t^\varphi(x))=(\tau_t\otimes \sigma_t^\varphi)\Delta(x)$,

\item[$ii)$] $\Delta(\tau_t(x))=(\tau_t\otimes\tau_t)\Delta(x)$,

\item[$iii)$] $\Delta(R(x))=(R\otimes R)\Delta'(x)$,
\end{enumerate}
where $\Delta'$ is obtained from~$\Delta$ by composing it with the flip
map on $M\otimes M$.
\end{Theorem}

\begin{proof}
As we saw already, the f\/irst formula is an immediate consequence of the second formula in Proposition~\ref{1.20} because
$\Delta(x)=W^*(1\otimes x)W$ for all $x\in M$.
The formulas ii) and iii) follow in a~straightforward way from the def\/initions and the formulas in
Proposition~\ref{1.20} and again $\Delta(x)=W^*(1\otimes x)W$, combined also with $(\Delta\otimes\iota)W=W_{13}W_{23}$.
For the third formula, take adjoints to get $(\Delta\otimes\iota)W*=W^*_{23}W^*_{13}$ and apply the f\/lip map on the
f\/irst two factors to get $(\Delta'\otimes\iota)W^*=W^*_{13}W^*_{23}$.
Then iii) follows from $W^*=(I\otimes J)W(I\otimes J)$ and the def\/inition of~$\tau_t$.
In fact, formula~ii) can also be obtained from~i) using coassociativity and the density~ii) in Proposition~\ref{1.21}.
\end{proof}

In order to complete the proof of the uniqueness of the left Haar weights, we also need the following result.

\begin{Proposition}\label{2.4}
If $x\in M$ and $\Delta(x)=1\otimes x$ or $\Delta(x)=x\otimes 1$, then~$x$ must be a~scalar multiple of~$1$.
\end{Proposition}

\begin{proof}
We will assume that $\Delta(x)=x\otimes 1$ and prove that~$x$ is a~scalar multiple of~$1$.
The other property (which is the one we really need) will follow by symmetry.

So assume that $x\in M$ and that $\Delta(x)=x\otimes 1$.
Take a~function~$f$ on ${\mathbb R}$ of the type $f(t)=\exp(-p(t-q)^2)$ with $p,q \in {\mathbb R}$ and def\/ine
\begin{gather*}
y=\int f(t) \sigma_t^\varphi(x)\, dt
\qquad
\text{and}
\qquad
z=\int f(t) \tau_t(x)\, dt.
\end{gather*}
Because $\Delta(x)=x\otimes 1$ and $\Delta(\sigma_t^\varphi(x))
=(\tau_t\otimes\sigma_t^\varphi)\Delta(x)=\tau_t(x)\otimes 1$, we get $\Delta(y)=z\otimes 1$.
Take any $w\in \mathcal M_\varphi$ ($=\mathcal N_\varphi^*\mathcal N_\varphi$).
Then, because~$y$ is analytic with respect to $(\sigma_t^\varphi)$, we have also $yw\in \mathcal M_\varphi$ (see e.g.\
Section~2 of Chapter~VIII in~\cite{T3}).
Therefore we can apply~$\varphi$ to the second leg of the equation $\Delta(yw)=(z\otimes 1)\Delta(w)$ and use left
invariance of~$\varphi$ to get $\varphi(yw)1=\varphi(w)z$.
Because this holds for all $w\in \mathcal M_\varphi$, it follows from the faithfulness of~$\varphi$ that~$y$ and~$z$ are
scalar multiples of the identity.
Because this is true for all such functions~$f$, this can only happen when~$x$ itself is a~scalar multiple of~$1$.
\end{proof}

Now we are ready to obtain the uniqueness of the Haar weights.

\begin{Theorem}
Any two left Haar weights on a~locally compact quantum group are equal (up to a~scalar).
Similarly for right Haar weights.
\end{Theorem}

\begin{proof}
Combining Proposition~\ref{2.2} where we have shown that $\Delta(u_t)=1\otimes u_t$ for all~$t$ and the previous result,
we f\/ind that $u_t$ is a~scalar multiple of $1$ for all~$t$.
This implies that the two weights $\varphi_1$ and $\varphi_2$ are proportional.

This result is not explicitly stated in Section~3 of Chapter~VIII in~\cite{T3}, but it follows easily as in the proof of
e.g.\ Corollary~3.6 in Chapter VIII of~\cite{T3}.
The result for right Haar weights follows by symmetry.
\end{proof}

Because of this result, we will in what follows use $(\sigma_t)$ for the modular automorphisms of a~left Haar weight (in
stead of $(\sigma_t^\varphi)$ as we did before).
We will use $(\sigma_t')$ for the modular automorphisms of a~right Haar weight.

From the fact that~$R$ f\/lips the coproduct, it follows that $\varphi\circ R$ is right invariant when~$\varphi$ is left
invariant.
From now on, we will assume that~$\varphi$ is a~{\it fixed left Haar weight} and we will take {\it the right Haar
weight}~$\psi$ to be this composition $\varphi\circ R$.

Now, here is another couple of formulas.

\begin{Theorem}\label{2.6}
For all $x\in M$ and $t\in {\mathbb R}$ we have
\begin{enumerate}\itemsep=0pt
\item[$i)$] $R(\sigma_t(x))=\sigma_{-t}'(R(x))$,

\item[$ii)$] $\Delta(\sigma_t'(x))=(\sigma_t' \otimes \tau_{-t})\Delta(x)$.
\end{enumerate}
\end{Theorem}

The f\/irst property is a~consequence of $\psi=\varphi\circ R$ and the fact that~$R$ is an anti-homomor\-phism.
It can be shown using e.g.\ the K.M.S.\ property of weights w.r.t.\ the modular automorphisms.
Having i), clearly ii) will follow from i) and iii) in Theorem~\ref{2.3}.
We could have obtained the second formula of this theorem also by symmetry, but then we would not be sure that we had
the same scaling group.

From the uniqueness of the Haar weights, we get easily that the scaling automorphisms $(\tau_t)$ leave the Haar weights
relatively invariant.

\begin{Theorem}\label{2.7}
There exists a~strictly positive number~$\nu$ so that $\varphi\circ\tau_t=\nu^{-t}\varphi$ and
$\psi\circ\tau_t=\nu^{-t}\psi$ for all $t\in {\mathbb R}$.
\end{Theorem}

\begin{proof}
It follows from $\Delta\circ\tau_t=(\tau_t\otimes\tau_t)\circ\Delta$ that $\varphi\circ\tau_t$ is also left invariant
and so, by uniqueness, it is a~scalar multiple of~$\varphi$.
As this is true for all $\tau_t$, we get a~strictly positive number~$\nu$ such that $\varphi\circ\tau_t=\nu^{-t}\varphi$
for all~$t$.
Composing with~$R$ and using that~$R$ and $\tau_t$ commute, we get also $\psi\circ\tau_t=\nu^{-t}\psi$ for all~$t$.
\end{proof}

From the fact that~$\varphi$ is relatively invariant under $\tau_t$, it follows that $\tau_t$ commutes with all the
modular automorphisms $(\sigma_s)_{s\in {\mathbb R}}$ of the left Haar weight.
This can be seen in dif\/ferent ways.
One can e.g.\ use the one-parameter group of unitaries $(v_t)$ def\/ined on $\mathcal H_\varphi$~by
\begin{gather*}
v_t\Lambda_\varphi(x)=\nu^{\frac12 t}\Lambda_\varphi(\tau_t(x))
\end{gather*}
and the fact that these unitaries commute with the map~$T$ (as def\/ined in Notation~\ref{1.19}).
The result then follows from the uniqueness of the polar decomposition of~$T$.
Similarly, $\tau_t$ will commute with all the modular automorphisms $(\sigma_s')_{s\in {\mathbb R}}$ of the right Haar
weight.
It is not so hard to get that then also the modular automorphisms of the left Haar weight commute with the modular
automorphisms of the right Haar weight.
Indeed, f\/ix $s,t\in {\mathbb R}$ and denote $\gamma=\sigma_t\tau_{-t}$ and $\gamma'=\sigma_s'\tau_s$.
From the formulas i) and ii) in Theorem~\ref{2.3} and ii) in Theorem~\ref{2.6} we get
\begin{gather*}
\Delta(\gamma(x)) =(\iota\otimes\gamma)\Delta(x),
\qquad
\Delta(\gamma'(x)) =(\gamma'\otimes \iota)\Delta(x)
\end{gather*}
for all~$x$.
It follows that $\Delta(\gamma\gamma'(x))=\Delta(\gamma'\gamma(x))$ for all~$x$.
This will imply $\gamma\gamma'=\gamma'\gamma$.
And as~$\tau$ and~$\sigma$ commute, as well as~$\tau$ and $\sigma'$, we can conclude from this that also~$\sigma$ and
$\sigma'$ will commute.
So we get the following result.

\begin{Theorem}\label{2.8}
All the automorphism groups~$\sigma$, $\sigma'$ and~$\tau$ mutually commute.
\end{Theorem}

Now we will show that~$\psi$ is relatively invariant w.r.t.\ the modular automorphism group~$\sigma$ and that~$\varphi$
is relatively invariant w.r.t.\ the modular automorphism group $\sigma'$, with (essentially) the same scaling
factor~$\nu$.
Observe that also this result would imply that~$\sigma$ and $\sigma'$ commute.
However, we will use that~$\sigma$ and $\sigma'$ commute to prove the following theorem.

\begin{Theorem}\label{2.9}
We have $\psi\circ\sigma_t=\nu^{-t}\psi$ and $\varphi\circ\sigma'_t=\nu^t \varphi$ for all $t\in {\mathbb R}$.
\end{Theorem}

\begin{proof}
Take $x\in \mathcal M_\psi$ and consider the formula $\Delta(\sigma_t(x))=(\tau_t\otimes\sigma_t)\Delta(x)$
(Theorem~\ref{2.3}.i).
If we also have that $\sigma_t(x)\in \mathcal M_\psi$, we can use $\psi\circ\tau_t=\nu^{-t}\psi$ (Theorem~\ref{2.7}),
and we f\/ind that $\psi(\sigma_t(x))=\nu^{-t}\psi(x)$ by invariance.
Because the weights $\psi\circ\sigma_t$ and $\nu^{-t}\psi$ have the same modular automorphism group (because~$\sigma$
and $\sigma'$ commute), it will follow (from Proposition~3.16 of Chapter~VIII in~\cite{T3}) that these weights are the
same if we can show that the $*$-subalgebra $\mathcal M_\psi \cap \sigma_{-t}(\mathcal M_\psi)$ is dense.
This is what we will do now.

Take any $x\in \mathcal N_\psi$ and take a~function~$f$ on ${\mathbb R}$ as in the proof of Proposition~\ref{2.4} to def\/ine
$y=\int f(s)\sigma_s'(x) ds$.
We know that~$y$ is still in $\mathcal N_\psi$ and analytic with respect to $\sigma'$.
Because~$\sigma$ commutes with $\sigma'$, also $\sigma_{-t}(y)$ will be analytic and so, as before, $\mathcal N_\psi
\sigma_{-t}(y)\subseteq \mathcal N_\psi$.
On the other hand, $\mathcal N_\psi \sigma_{-t}(y) \subseteq \sigma_{-t}(\mathcal N_\psi)$ because $y\in \mathcal
N_\psi$ and $\mathcal N_\psi$ is a~left ideal.
So we have produced elements in the intersection of $\mathcal N_\psi$ and $\sigma_{-t}(\mathcal N_\psi)$.
It is not hard to conclude that this intersection is dense, as well as the intersection $\mathcal M_\psi\cap\sigma_{-t}(\mathcal M_\psi)$.
This completes the proof of the f\/irst part of the theorem.

The second formula is obtained from the f\/irst one by using i) of Theorem~\ref{2.6}.
\end{proof}

Remark that the proof of this result is dif\/ferent from the original one.
In~\cite{K-V2}, f\/irst a~stronger form of right invariance for~$\psi$ is needed.

Now we can add one more relation of the type proven in Theorems~\ref{2.3} and~\ref{2.6}.

\begin{Theorem}\label{2.10}
For all $x\in M$ we have $\Delta(\tau_t(x))=(\sigma_t\otimes\sigma'_{-t})\Delta(x)$ for all $t\in {\mathbb R}$.
\end{Theorem}

\begin{proof}
Because of the relative invariance of~$\varphi$ under both $\tau_t$ and $\sigma'_t$, we have one-parameter groups of
unitaries $(v_t)$ and $(w_t)$ on $\mathcal H_\varphi$ satisfying
\begin{gather*}
v_t\Lambda_\varphi(x) =\nu^{\frac12 t}\Lambda_\varphi(\tau_t(x)),
\qquad
w_t\Lambda_\varphi(x) =\nu^{-\frac12 t}\Lambda_\varphi(\sigma'_t(x))
\end{gather*}
for all $x\in \mathcal N_\varphi$ and all~$t$.
When~$W$ is the left regular representation as before, it follows from the equations
$\Delta\circ\tau_t=(\tau_t\otimes\tau_t)\Delta$ and $\Delta\circ\sigma'_t=(\sigma'_t\otimes\tau_{-t})\Delta$ that
\begin{gather*}
(1\otimes v_t^*)W^*(1\otimes v_t)  = (\tau_t\otimes\iota)W^*,
\qquad
(1\otimes v_t)W^*(1\otimes w_t)  = (\sigma'_t\otimes\iota)W^*.
\end{gather*}
If we apply $R\otimes J(\,\cdot\,)^*J$ to the second equation, we get because~$J$ commutes with $w_t$ and $v_t$ that
\begin{gather*}
(1\otimes w_t^*)W^*(1\otimes v_t^*) = (\sigma_{-t}\otimes\iota)W^*.
\end{gather*}
If we combine all of this with the formula $(\Delta\otimes\iota)W=W_{13}W_{23}$ as we did before, we can conclude that
$\Delta\circ\tau_t=(\sigma_t\otimes\sigma'_{-t})\Delta$ on the left leg of~$W$.
As this leg is dense, we get the result.
\end{proof}

The unitary groups $(v_t)$ and $(w_t)$, def\/ined in the proof of this theorem, will play an important role further.
This will be seen in Section~\ref{Section5} where we will recall the def\/initions.

From the fact that~$\psi$ is relatively invariant under the modular automorphisms of~$\varphi$, we also get the
following important result.

\begin{Theorem}\label{2.11}
There exists a~unique, non-singular, positive self-adjoint operator~$\delta$, affiliated with~$M$ such that
$\psi=\varphi\big(\delta^{\frac12}\,\cdot\,\delta^{\frac12}\big)$.
This operator satisfies $\sigma_t(\delta)=\nu^t\delta$ and $\sigma'_t(\delta)=\nu^t\delta$.
It is invariant under the automorphisms $(\tau_t)$ and $R(\delta)=\delta^{-1}$.
We also have the relation $\sigma'_t(x)=\delta^{it}\sigma_t(x)\delta^{-it}$.
\end{Theorem}

First remark that we have formulated the above results in terms of the unbounded operator~$\delta$.
It is quite obvious what is meant by this for all the formulas in the formulation, except for the f\/irst one.
This will have to be interpreted with the use of the Connes' cocycle Radon--Nikodym derivative as we will see in the
proof.
We refer to Corollary~3.6 in Section~VIII of~\cite{T3} and to~\cite{V1} for more details.
In fact, when possible, we will rather consider the equivalent formulas in terms of the unitary operators
$(\delta^{it})$ in order to avoid this kind of (technical) dif\/f\/iculties.

\begin{proof}
Consider the weights~$\varphi$ and $\psi=\varphi\circ R$ and consider the cocycle Radon--Nikodym derivative.
We will write $u_t=(\text{D}\psi:\text{D}\varphi)_t$ for all~$t$.
For any automorphism~$\alpha$ of~$M$, we have
\begin{gather*}
(\text{D}\psi\circ\alpha:\text{D}\varphi\circ\alpha)_t=\alpha^{-1}(\text{D}\psi:\text{D}\varphi)_t
\end{gather*}
and if we apply this with $\alpha=\sigma_s$ and if we use that $\varphi\circ\sigma_s=\varphi$ and
$\psi\circ\sigma_s=\nu^{-s}\psi$, we get $\sigma_s(u_t)=\nu^{ist}u_t$ for all~$t$ and all~$s$.
Because $u_{s+t}=u_s\sigma_s(u_t)$ we easily calculate that $(\nu^{-\frac12 it^2}u_t)$ is a~one-parameter group of
unitaries in~$M$.
Therefore, there exists a~non-singular positive self-adjoint operator~$\delta$, af\/f\/iliated with~$M$, such that
$u_t=\nu^{\frac12 it^2}\delta^{it}$ for all~$t$.
This gives us the element~$\delta$ such that, at least formally,
$\psi=\varphi(\delta^{\frac12}\,\cdot\,\delta^{\frac12})$.
Takesaki (Chapter VIII in~\cite{T3})
treats this situation in the case $\nu=1$ whereas Vaes~\cite{V1} considers the general case.

Because~$\psi$ and~$\varphi$ are scaled with the same factor by the automorphisms $\tau_s$ we must have that $u_t$ and
hence also~$\delta$ is invariant under $\tau_s$.
Because $\psi=\varphi\circ R$ and $\varphi\circ R=\psi$ we will get that $R(u_t)=u_{-t}$ and so $R(\delta)=\delta^{-1}$.
Because $\sigma_t(u_s)=\nu^{ist}u_s$ we obtain $\sigma_t(\delta)=\nu^t\delta$ and f\/inally, if we apply~$R$ to this last
formula, we get the other one $\sigma'_t(\delta)=\nu^t\delta$.

The last statement is standard.
\end{proof}

One also has the formula $\Delta(\delta)=\delta \otimes \delta$ but this is not so easy to get.
One possible way to prove this formula can be found in the original work~\cite{K-V2}.
Another proof is found in~\cite{M-N-W}.
We will still give an other argument but for this, we need some more results and so we postpone the proof of this
formula (see Remark~\ref{4.16}).

It is possible to realize the GNS-representation of~$\psi$ in the Hilbert space $\mathcal H_\varphi$.
One can show that any element $x\in M$ with the property that $x\delta^{\frac12}$ is bounded and belongs to $\mathcal
N_\varphi$ is an element of $\mathcal N_\psi$.
In that case $\Lambda_\psi(x)=\Lambda_\varphi(x\delta^{\frac12})$.
Moreover, it is not hard to prove that there are enough elements like that and that $\Lambda_\psi$ is completely
determined in this way.
All this is proven by left Hilbert algebra techniques.
Details can be found e.g.\ in~\cite{V1}, see also~\cite{VD8}.

\section[The dual $(\widehat M,\widehat\Delta)$]{The dual $\boldsymbol{(\widehat M,\widehat\Delta)}$}\label{Section4}

In this section, {\it we start with a~locally compact quantum group} $(M,\Delta)$ as in Def\/inition~\ref{2.1} and we will
construct the dual $(\widehat M,\widehat\Delta)$.
We will also show that repeating the procedure gives the original pair $(M,\Delta)$.

We take a~left Haar weight~$\varphi$ and the right Haar weight~$\psi$ obtained from it by composing with the unitary
antipode~$R$ (see the~remark before Theorem~\ref{2.6} in the previous section).
We will, as explained also at the end of the previous section, identify $\mathcal H_\psi$ with $\mathcal H_\varphi$~by
def\/ining $\Lambda_\psi$ from $\mathcal N_\psi$ to $\mathcal H_\varphi$~by
$\Lambda_\psi(x)=\Lambda_\varphi(x\delta^{\frac12})$ for the appropriate elements~$x$.
Observe again that this is compatible with the actions of~$M$ on both spaces.
Moreover, because we have uniqueness of the Haar weights, it is no longer necessary to use the subscript~$\varphi$ when
considering the Hilbert space.
So, from now on, we will simply use $\mathcal H$ for $\mathcal H_\varphi$.
In other words, we will be working solely within the Hilbert space $\mathcal H_\varphi$ and denote it by $\mathcal H$.

We consider the left regular representation~$W$ of $(M,\Delta)$ associated with the left Haar weight as in Proposition~\ref{1.11}.
With the above conventions, we get that~$W$ acts on $\mathcal H\otimes \mathcal H$ and that $W\in M \otimes \mathcal B(\mathcal H)$.
Then also~$W$ satisf\/ies the pentagon equation $W_{12}W_{13}W_{23}=W_{23}W_{12}$.

We will make use of the antipode~$S$ and its polar decomposition $S=R\tau_{-\frac{i}{2}}$.
We will also use that the space of elements of the form $(\iota\otimes\omega)W$ with $\omega\in \mathcal B(\mathcal
H)_*$ is invariant under the scaling group $(\tau_t)$, that it is a~core for~$S$ and that
$S((\iota\otimes\omega)W)=(\iota\otimes\omega)W^*$.
We will refer to this property when we write formally $(S\otimes\iota)W=W^*$.
For details see Proposition~\ref{1.24}.

\medskip

{\bf The dual von Neumann algebra and the dual coproduct.}
First the underlying von Neumann algebra $\widehat M$ is def\/ined.
\begin{Definition}\label{3.1}
Let $\widehat M$ be the~$\sigma$-weak closure of the subspace $\{(\omega\otimes\iota)W \,|\, \omega\in M_*\}$ of
$\mathcal B(\mathcal H)$.
\end{Definition}

With this def\/inition, the dual von Neumann algebra $\widehat M$ also acts on the space $\mathcal H$ (which is~$\mathcal
H_\varphi$).
In some sense, this means that from the very beginning, we identify the `spaces' $L^2(G)$ and $L^2(\widehat G)$, thus
not using the `Fourier transform' explicitly.
We come back to this remark later in this section when we construct the dual left Haar weight.

It follows from the fact that~$W$ is a~multiplicative unitary, that the subspace~$C$ of $\mathcal B(\mathcal H)$,
def\/ined as $\{(\omega\otimes\iota)W \,|\, \omega\in M_*\}$, is a~subalgebra of $\mathcal B(\mathcal H)$.
In order to have that its closure is a~von Neumann algebra, we will use the following result (which will also be needed
further in this section).

\begin{Lemma}\label{3.2}
Consider the one-parameter group of automorphisms $(\tau_t)$ on~$M$ and its dual on~$M_*$.
If $\omega\in M_*$ is analytic and if $\omega_1$ is defined as $\overline{\omega}\circ\tau_{-\frac{i}{2}}\circ R$, then
\begin{gather*}
((\omega\otimes\iota)W)^*=(\omega_1\otimes\iota)W.
\end{gather*}
Conversely, if $\omega,\omega_1\in M_*$ and satisfy the above equation, then
$\omega_1(x)=\omega(S(x)^*)^-=\overline\omega(S(x))$ for all $x\in \mathcal D(S)$.
\end{Lemma}

\begin{proof}
Recall that formally $(S\otimes\iota)W=W^*$ and that $S=R\tau_{-\frac{i}{2}}$.
Now, it is not so hard to show rigorously that the f\/irst statement of the lemma is correct.

Conversely, let $\omega,\omega_1\in M_*$ and assume that $((\omega\otimes\iota)W)^*=(\omega_1\otimes\iota)W$.
If $x=(\iota\otimes\rho)W$ with $\rho\in \mathcal B(\mathcal H)_*$, then
\begin{gather*}
\omega_1(x) = \rho((\omega_1\otimes\iota)W) = \overline\rho((\omega\otimes\iota)W)^-
 = \omega((\iota\otimes\overline\rho)W)^- = \omega(S(x)^*)^-.
\end{gather*}
Now, because the `left leg' of~$W$ is a~core for~$S$, we will also get the desired formula for all $x\in \mathcal D(S)$.
\end{proof}

Now the following (standard) result can be shown.

\begin{Proposition}\label{3.3}
$\widehat M$ is a~von Neumann algebra and $W\in M\otimes\widehat M$.
\end{Proposition}

\begin{proof}
Again consider the subspace~$C$ def\/ined by $\{(\omega\otimes\iota)W \,|\, \omega\in M_*\}$ (without taking the closure).
By the f\/irst statement of the lemma and using that such analytic elements are dense, we get that the closure $\widehat M$ is self-adjoint.
It acts non-degenerately on $\mathcal H$ because~$W$ is a~unitary.
Therefore we have that $\widehat M$ will be a~von~Neumann algebra on $\mathcal H$.

The second statement follows from the commutation theorem for tensor products of von Neumann algebras.
\end{proof}

We see from the above that also the norm closure of the space~$C$ is a~$C^*$-algebra, acting non-degenerately on the Hilbert space $\mathcal H$.
It is denoted by $\widehat A$ and it will be considered further in  Appendix~\ref{appendixA} where we treat the relation between the
von~Neumann
algebra approach and the $C^*$-algebra approach.

Remark that it is not immediate that these subalgebras $\widehat A$ and $\widehat M$ are self-adjoint.
Usually either regularity (cf.~\cite{B-S}) or manageability (cf.~\cite{W2}) of the multiplicative unitary~$W$ is used.
Here we use a~result which is closely related to manageability, but (in some sense) also weaker than manageability.
As explained already in the introduction, we avoid the use of the concept of manageability.

Finally we also want to mention the following.
Because $\Delta(x)=W(1\otimes x)W^*$ for $x\in M$, it follows from Proposition~\ref{2.4} that $M\cap \widehat M'={\mathbb C} 1$.
Indeed, if $x\in M\cap \widehat M'$ then $W(1\otimes x)=(1\otimes x)W$ because $x\in \widehat M'$ and so
$\Delta(x)=1\otimes x$.
By Proposition~\ref{2.4} we get $x\in {\mathbb C} 1$.

Now we proceed to def\/ine the coproduct $\widehat \Delta$ on $\widehat M$.
We f\/irst formulate the following result.

\begin{Lemma}\label{3.4}

We have $W(y\otimes 1)W^*\in \widehat M\otimes\widehat M$ for all $y\in \widehat M$.
\end{Lemma}

This result is standard and easy to prove.
It follows from the def\/inition of $\widehat M$, the formula $W_{23}W_{12}W_{23}^*=W_{12}W_{13}$ and the fact that $W\in M\otimes \widehat M$.

Now it is easy to def\/ine the coproduct $\widehat \Delta$.
As is common in the theory of locally compact quantum groups, we use the f\/lip in the following def\/inition.

\begin{Definition}\label{3.5}
We let $\widehat\Delta(y)=\chi(W(y\otimes 1)W^*)$ for $y\in \widehat M$ where~$\chi$ is the f\/lip on the tensor product
$\widehat M\otimes \widehat M$.
\end{Definition}

It is clear, using the previous lemma, that $\widehat \Delta$ is a~unital and normal $*$-homomorphism from~$\widehat M$
to $\widehat M\otimes \widehat M$.
The coassociativity follows in the standard way from the pentagon equation.

\medskip

{\bf Construction of the left Haar weight $\boldsymbol{\widehat\varphi}$ on $\boldsymbol{(\widehat M,\widehat\Delta)}$.}
The next step is the construction of the Haar weights on the dual $(\widehat M,\widehat \Delta)$.
We will f\/irst construct the left Haar weight $\widehat\varphi$ and later discuss the existence of the right Haar weight $\widehat\psi$.
We will follow the (more or less) standard procedure.

We begin with the construction of the map $\widehat \Lambda$ (which later will be closed and give the map
$\Lambda_{\widehat\varphi}$ associated with the dual left Haar weight $\widehat\varphi$).
The def\/inition looks somewhat strange although it is well-known in the theory of Kac algebras.
Let us give here a~short motivation using the algebraic theory of multiplier Hopf algebras with integrals~\cite{VD4}.

In~\cite{VD4} the Fourier transform $\widehat a$ of an element~$a$ is def\/ined as the linear functional
$\omega=\varphi(\,\cdot\, a)$.
Now it is well-known that the multiplicative unitary~$W$ is essentially the duality (see e.g.~\cite{VD1} for the Hopf
algebra case and~\cite{D-VD} for the case of algebraic quantum groups).
So, formally, we have $\widehat a=(\omega\otimes\iota)W$ when $\omega=\varphi(\,\cdot\, a)$.
As we remarked already before, the `spaces'~$L^2(G)$ and~$L^2(\widehat G)$ are identif\/ied which means (again formally)
that we want $\Lambda_{\widehat\varphi}(\widehat a)=\Lambda_\varphi(a)$.
This formula is rewritten as
\begin{gather*}
\langle\Lambda_{\widehat\varphi}(\widehat a),\Lambda_\varphi(x)\rangle =
\langle\Lambda_\varphi(a),\Lambda_\varphi(x)\rangle =\varphi(x^*a)=\omega(x^*),
\end{gather*}
whenever $x\in \mathcal N_\varphi$.

Therefore the following def\/inition is not a~surprise.

\begin{Definition}\label{3.6}
Def\/ine a~subspace $\widehat{\mathcal N}$ of $\widehat M$ and a~linear map $\widehat\Lambda:\widehat{\mathcal N} \to
\mathcal H$ as follows.
We say that an element~$y$ of $\widehat M$ belongs to $\widehat{\mathcal N}$ if there exists a~linear functional
$\omega\in M_*$ such that $y=(\omega\otimes\iota)W$ and a~vector $\xi\in \mathcal H$ such that
$\langle\xi,\Lambda_\varphi(x)\rangle=\omega(x^*)$ for all $x\in \mathcal N_\varphi$.
Then we set $\widehat\Lambda(y)=\xi$.
\end{Definition}

As before, $\Lambda_\varphi$ is the canonical map from $\mathcal N_\varphi$ to $\mathcal H$.
Remark that~$\omega$ is uniquely determined by~$y$ because the left leg of~$W$ is dense in~$M$ (cf.\ Proposition~\ref{1.21})
and that $\xi$ is completely determined by~$\omega$ because $\Lambda_\varphi(\mathcal N_\varphi)$ is dense in $\mathcal H$.
It is also clear that $\widehat {\mathcal N}$ is a~subspace and that~$\widehat \Lambda$ is linear and injective.

In the next lemma, we will show that there are enough elements in this space $\widehat {\mathcal N}$.
We will need {\it right bounded} vectors in $\mathcal H$.
Recall that these are elements $\eta\in\mathcal H$ such that there exist a~bounded operator, denoted by $\pi'(\eta)$,
satisfying $\pi'(\eta)\Lambda_\varphi(x)=x\eta$ for all $x\in \mathcal N_\varphi$.
Such vectors form a~dense subspace and the space of operators $\pi'(\eta)$ with~$\eta$ right bounded, is dense in the
commutant $M'$ of~$M$ (see e.g.\ Chapter VI in~\cite{T3}).

\begin{Lemma}\label{3.7}
Let $\xi,\eta\in \mathcal H$ and assume that~$\eta$ is right bounded.
Let $\omega=\langle \,\cdot\,\xi,\eta\rangle$ and $y=(\omega\otimes\iota)W$.
Then $y\in \widehat{\mathcal N}$ and $\widehat\Lambda(y)=\pi'(\eta)^*\xi$.
In particular, $\widehat{\mathcal N}$ is~$\sigma$-weakly dense in $\widehat M$ and the space
$\widehat\Lambda(\widehat{\mathcal N})$ is dense in $\mathcal H$.
\end{Lemma}

\begin{proof}
If $\xi,\eta\in \mathcal H$ and if~$\eta$ is right bounded, we have for all $x\in \mathcal N_\varphi$ that
\begin{gather*}
\omega(x^*)=\langle x^*\xi,\eta\rangle=\langle\xi,x\eta\rangle=\langle \xi,\pi'(\eta)\Lambda_\varphi(x)\rangle=\langle
\pi'(\eta)^*\xi,\Lambda_\varphi(x)\rangle.
\end{gather*}
So, if $y=(\omega\otimes \iota)W$, then $y\in \widehat{\mathcal N}$ and $\widehat\Lambda(y)=\pi'(\eta)^*\xi$.

Because the space of right bounded vectors is dense in $\mathcal H$, we see that $\widehat{\mathcal N}$ is dense in $\widehat M$.
And because the space of operators $\pi'(\eta)$ with~$\eta$ right bounded is dense in $M'$, we get that also the space
$\widehat\Lambda(\widehat{\mathcal N})$ is dense in $\mathcal H$.
\end{proof}

We also have the following.

\begin{Lemma}\label{3.8}
Let $\omega,\omega_1\in M_*$ and $y=(\omega\otimes\iota)W$ and $y_1=(\omega_1\otimes\iota)W$.
If $y\in \widehat{\mathcal N}$ then also $y_1y\in \widehat{\mathcal N}$ and $\widehat\Lambda(y_1y)=y_1\widehat \Lambda(y)$.
\end{Lemma}

\begin{proof}
For any $x\in\mathcal N_\varphi$ we have
\begin{gather*}
\langle y_1 \widehat\Lambda(y), \Lambda_\varphi(x) \rangle = \langle \widehat\Lambda(y), y_1^*\Lambda_\varphi(x)\rangle
 = \langle \widehat\Lambda(y), ((\overline{\omega_1}\otimes\iota)(W^*))\Lambda_\varphi(x) \rangle
\\
\phantom{\langle y_1 \widehat\Lambda(y), \Lambda_\varphi(x) \rangle}
 = \langle \widehat\Lambda(y), \Lambda_\varphi((\overline{\omega_1}\otimes\iota)\Delta(x)) \rangle
 = \omega(((\overline{\omega_1}\otimes\iota)\Delta(x))^*)
\\
\phantom{\langle y_1 \widehat\Lambda(y), \Lambda_\varphi(x) \rangle}
 = \omega((\omega_1\otimes\iota)\Delta(x^*))
 = (\omega_1\omega)(x^*),
\end{gather*}
where $\omega_1\omega$ is def\/ined as usual by $(\omega_1\omega)(a)=(\omega_1\otimes\omega)\Delta(a)$ for all $a\in M$.
Remark that we have used the def\/inition of $W^*$ as given in Proposition~\ref{1.11}.
Finally it is easy to see that $y_1y=((\omega_1 \omega)\otimes\iota)W$ using the pentagon equation and the formula
$\Delta(a)=W^*(1\otimes a)W$ for $a\in M$.
\end{proof}

The left Haar weight $\widehat\varphi$ will be obtained by constructing a~left Hilbert algebra.
Def\/inition~\ref{3.6} and Lemmas~\ref{3.7} and~\ref{3.8} clearly provide the f\/irst steps.
We need one more lemma before we can come to the main part of the construction.

\begin{Lemma}
Consider the set $\widehat{\mathcal N}_0$ of elements $y\in \widehat{\mathcal N}$ so that $y^*$ has the form
$(\omega_1\otimes\iota)W$ for some $\omega_1\in M_*$.
Then $\widehat{\mathcal N}_0$ is still~$\sigma$-weakly dense in $\widehat M$ and also $\widehat\Lambda(\widehat{\mathcal
N}_0)$ is still dense in $\mathcal H$.
\end{Lemma}

\begin{proof}
Take $\omega=\langle \,\cdot\,\xi,\eta\rangle$ with $\xi,\eta\in \mathcal H$ and assume that~$\eta$ is right bounded.
We know from Lemma~\ref{3.7} that~$y$, def\/ined as $(\omega\otimes\iota)W$, is in $\widehat{\mathcal N}$ and that
$\widehat\Lambda(y)=\pi'(\eta)^*\xi$.
Now, because of Lemma~\ref{3.2}, we need such elements with~$\omega$ analytic with respect to $(\tau_t)$.

To construct such elements, def\/ine a~one-parameter group of unitaries $(v_t)$ on $\mathcal H$ (as in the proof of
Theorem~\ref{2.10}) by $v_t\Lambda_\varphi(x)=\nu^{\frac12 t}\Lambda_\varphi(\tau_t(x))$ when $x\in \mathcal N_\varphi$.
It is clear that $\tau_t(x)=v_t x v_t^*$ for all $x\in M$.
It will be possible to take the vector $\xi$ above so that it is analytic with respect to $(v_t)$.
And because every $v_t$ will map right bounded elements to right bounded elements, we will also be able to take~$\eta$
analytic and still right bounded.
Then~$\omega$ will be analytic and of the required form.
This will give us the result of the lemma.
\end{proof}

Now we come to the main step in the construction of the left Hilbert algebra we need to get the dual weight
$\overline\varphi$.

\begin{Proposition}
Let $\mathfrak A= \widehat\Lambda(\widehat{\mathcal N}\cap\widehat{\mathcal N}^*)$.
We can equip $\mathfrak A$ with the $*$-algebra structure inherited from $\widehat{\mathcal N}\cap\widehat{\mathcal N}^*$.
If we denote~$y$ by $\pi(\xi)$ when $y\in \widehat{\mathcal N}\cap\widehat{\mathcal N}^*$ and $\xi=\widehat\Lambda(y)$,
then we have
\begin{enumerate}\itemsep=0pt
\item[$i)$] both $\mathfrak A$ and $\mathfrak A^2$ are dense in $\mathcal H$,

\item[$ii)$] $\pi(\xi)$ is a~bounded operator for all $\xi\in \mathfrak A$,

\item[$iii)$] $\pi$ is a~$*$-representation of $\mathfrak A$.
\end{enumerate}
\end{Proposition}

\begin{proof}
It follows from the def\/inition of $\widehat{\mathcal N}$ and Lemma~\ref{3.8} that $\widehat{\mathcal N}$ is a~subalgebra
of $\widehat M$ and so $\widehat{\mathcal N}\cap\widehat{\mathcal N}^*$ is a~$*$-subalgebra of $\widehat M$.
Because $\widehat\Lambda$ is injective, we can equip $\mathfrak A$ with the $*$-algebra structure of $\widehat{\mathcal N}\cap\widehat{\mathcal N}^*$.

To prove i), we use the previous lemma.
Indeed, if we take $y$, $y_1$ in the set $\widehat{\mathcal N}_0$, we see that both $y^*y_1$ and $y_1^*y$ will be in
$\widehat{\mathcal N}$ and this will provide us with enough elements in $\widehat{\mathcal N}\cap\widehat{\mathcal N}^*$.
And because this set is dense in $\widehat M$ and also $\widehat\Lambda(\widehat{\mathcal N}\cap\widehat{\mathcal N}^*)$
is dense, we get i).

Statements ii) and iii) are immediate consequences of the notations.
\end{proof}

There is now only one point missing for $\mathfrak A$ to be a~left Hilbert algebra.
It is also needed that the $*$-operation in $\mathfrak A$, usually denoted by $\xi\mapsto\xi^\sharp$, is preclosed.
This will be a~consequence of the following lemma.

\begin{Lemma}\label{3.11}
If $y\in \widehat{\mathcal N}\cap \widehat{\mathcal N}^*$ then
\begin{gather*}
\langle\widehat\Lambda(y^*),\Lambda_\varphi(a)\rangle= \langle\widehat\Lambda(y),\Lambda_\varphi(S(a^*))\rangle^-,
\end{gather*}
whenever $a\in \mathcal N_\varphi$, $a^*\in \mathcal D(S)$ and $S(a^*)\in \mathcal N_\varphi$.
\end{Lemma}

\begin{proof}
Take $y\in \widehat{\mathcal N}\cap\widehat{\mathcal N}^*$ and let $\omega$, $\omega_1$ be in $M_*$ so that
$y=(\omega\otimes\iota)W$ and $y^*=(\omega_1\otimes\iota)W$.
From Lemma~\ref{3.2} we know that $\omega_1(x)=\omega(S(x)^*)^-$ for all $x\in \mathcal D(S)$.
If now~$a$ is as in the formulation of the lemma, it will follow from the def\/inition of $\widehat\Lambda$ that
\begin{gather*}
\langle\widehat\Lambda(y^*),\Lambda_\varphi(a)\rangle = \omega_1(a^*)=\omega(S(a^*)^*)^-
 = \langle\widehat\Lambda(y),\Lambda_\varphi(S(a^*))\rangle^-
\end{gather*}
and the result will follow.
\end{proof}

So in order to prove that the involution in $\mathfrak A$ is preclosed, we just need to argue that there are enough
elements~$a$ as in the lemma.
This is the content of the next lemma.

\begin{Lemma}\label{3.12}
The set of elements $\Lambda_\varphi(a)$ with $a\in \mathcal N_\varphi$ such that also $a^*\in \mathcal D(S)$ and
$S(a^*)\in \mathcal N_\varphi$ is dense in $\mathcal H$.
\end{Lemma}

\begin{proof}
Take $a\in \mathcal N_\varphi$.
Formally we have
\begin{gather*}
S(a^*) =R \tau_{-\frac{i}{2}}(a^*)=R(\tau_{\frac{i}{2}}(a))^*
 =R(\tau_{\frac{i}{2}}(a)\delta^\frac{1}{2})^*\delta^\frac{1}{2}.
\end{gather*}
For such an element to be again in $\mathcal N_\varphi$, we need that $\tau_{\frac{i}{2}}(a)\delta^\frac12$ is
well-def\/ined and in $\mathcal N_\varphi$.
We need two results to obtain such elements.
First observe that~$\varphi$ is relatively invariant w.r.t.\ the automorphisms $(\tau_t)$ and so also $\mathcal N_\varphi$
is invariant and standard techniques allow to produce elements $a\in \mathcal N_\varphi$ that are analytic w.r.t.\ $(\tau_t)$.
Next we know that the elements $\delta^{is}$ are analytic w.r.t.\ the modular automorphisms $(\sigma_t)$ and so
$\mathcal N_\varphi \delta^{is}\subseteq \mathcal N_\varphi$ for all~$s$.
This will allow us to produce elements $a\in \mathcal N_\varphi$ such that $a\delta^\frac12$ is well-def\/ined and still in $\mathcal N_\varphi$.
The two techniques together will give us enough of the desired elements.
\end{proof}

We know from Section~\ref{Section2} that the operator~$K$ is essentially the map $\Lambda_\psi(x)\mapsto \Lambda_\psi(S(x)^*)$.
A~simple (formal) argument then shows that the map $\Lambda_\varphi(x)\mapsto \Lambda_\varphi(S(x^*))$ is essentially the operator $K^*$.
So we expect that the map $\widehat\Lambda(y)\mapsto \widehat\Lambda(y^*)$ will be nothing else but the operator~$K$.
It does not seem to be easy to prove this result exactly.
In the next section, we will f\/ind a~way around this problem by (in some sense) `redef\/ining' these maps (see Remark~\ref{4.10}).
We will then also take up again the argument that we sketch here in the proof of Lemma~\ref{3.12} (see the proof of Theorem~\ref{4.11}).

For the moment, there is no need to get this more precise result.
Indeed, if we combine all the previous results, we f\/ind that $\mathfrak A$ is a~left Hilbert algebra and this is what we need.
Then we can use the general procedure to construct a~faithful normal semi-f\/inite weight from a~left Hilbert algebra (see Chapter~VII in~\cite{T3})
and we will arrive at the following theorem.

\begin{Theorem}\label{3.13}
There exists a~normal faithful semi-finite weight $\widehat\varphi$ on $\widehat M$ such that the G.N.S.-representation
can be realized in $\mathcal H$, satisfying $\widehat{\mathcal N}\subseteq \mathcal N_{\widehat\varphi}$ and such that
the canonical map $\Lambda_{\widehat \varphi}$ is the closure of $\widehat \Lambda$ on $\widehat{\mathcal N}$.
\end{Theorem}

\begin{proof}
We have the left Hilbert algebra $\widehat\Lambda(\widehat{\mathcal N}\cap \widehat{\mathcal N}^*)$ sitting inside
$\mathcal H$.
The canonical weight~$\widehat\varphi$ associated to this left Hilbert algebra has the property that elements in
$\widehat{\mathcal N}\cap \widehat{\mathcal N}^*$ belong to~$\mathcal N_{\widehat\varphi}$ and that
$\widehat\varphi(y^*y)=\langle\widehat\Lambda(y),\widehat\Lambda(y)\rangle$ for such elements.
The rest follows from standard Hilbert algebra theory and the construction method of the associated weight.
\end{proof}

\medskip

{\bf The left invariance of $\boldsymbol{\widehat\varphi}$ and the associated regular representation.}
Now we need to show that $\widehat\varphi$ is left invariant on the pair $(\widehat M,\widehat \Delta)$.
When this is done, we can easily construct the right Haar weight on $(\widehat M,\widehat \Delta)$.
Indeed, as in the proof of Theorem~\ref{2.3}, we will have that $\widehat \Delta(\widehat R(y))=\chi(\widehat
R\otimes\widehat R)\widehat\Delta(y)$ whenever $y\in \widehat M$ where as before~$\chi$ denotes the f\/lip and where here
$\widehat R$ is def\/ined on $\widehat M$ by $\widehat R(y)=Jy^*J$.
Recall that~$J$ is the modular conjugation associated with~$\varphi$ as def\/ined in Notation~\ref{1.19}.
Therefore a~right Haar weight can be constructed from the left Haar weight on $\widehat M$ by composing it with this map $\widehat R$.
We will show later, in Section~\ref{Section5}, that the use of the notation $\widehat R$ is justif\/ied (see Proposition~\ref{4.12}).

This will give us that the pair $(\widehat M,\widehat \Delta)$ is again a~locally compact quantum group.
In order to show that repeating the procedure will bring us back to the original locally compact quantum group
$(M,\Delta)$, we will prove that the left regular representation of the dual is nothing else but the unitary $\Sigma
W^*\Sigma$ where as is common,~$\Sigma$ denotes the f\/lip operator on the tensor product $\mathcal H\otimes \mathcal H$.

We will prove this result f\/irst because it can be used to show that the dual left Haar weight~$\widehat \varphi$ is indeed left invariant.
In other words, the main result left to prove is the following.

\begin{Proposition}\label{3.14}
Define the unitary $\widehat W=\Sigma W^* \Sigma$ on $\mathcal H\otimes \mathcal H$.
Then $(\omega\otimes\iota)\widehat\Delta(y)\in \mathcal N_{\widehat\varphi}$ and
\begin{gather*}
\big((\omega\otimes\iota)\widehat
W^*\big)\Lambda_{\widehat\varphi}(y)=\Lambda_{\widehat\varphi}\big((\omega\otimes\iota)\widehat\Delta(y)\big),
\end{gather*}
whenever $y\in \mathcal N_{\widehat\varphi}$ and $\omega\in \mathcal B(\mathcal H)_*$.
\end{Proposition}

\begin{proof}
First take $y\in \widehat{\mathcal N}$ and let $\omega\in M_*$ be such that $y=(\omega\otimes\iota)W$.
Take any $\rho\in \widehat M_*$ and put $y_1=(\rho\otimes\iota)\widehat \Delta(y)$.
Since $\widehat\Delta(y)=\chi(W(y\otimes 1)W^*)$ (cf.\ Def\/inition~\ref{3.5}), it follows from a~straightforward
calculation (using the pentagon equation) that $y_1=(\omega(\,\cdot\,c)\otimes\iota)W$ where $c=(\iota\otimes\rho)W$.
Then, when $a\in \mathcal N_\varphi$, we have
\begin{gather*}
\omega(a^*c) = \langle\widehat\Lambda(y),c^*\Lambda_\varphi(a)\rangle
 = \langle c\widehat\Lambda(y),\Lambda_\varphi(a)\rangle
\end{gather*}
and it follows that $y_1\in\widehat{\mathcal N}$ and that $\widehat\Lambda(y_1)=c\widehat\Lambda(y)$.
This precisely means that
\begin{gather*}
\widehat\Lambda\big((\rho\otimes\iota)\widehat\Delta(y)\big) = ((\iota\otimes\rho)W)\widehat\Lambda(y)
 = \big((\rho\otimes\iota)\widehat W^*\big)\widehat\Lambda(y).
\end{gather*}
This proves the result when $y\in \widehat{\mathcal N}$.
The general case follows because $\Lambda_{\widehat \varphi}$ on $\mathcal N_{\widehat\varphi}$ is the closure of
$\widehat \Lambda$ on $\widehat{\mathcal N}$.
\end{proof}

Now it is not hard to prove left invariance of $\widehat\varphi$.

\begin{Proposition}\label{3.15}
The weight $\widehat\varphi$, as constructed in Proposition~{\rm \ref{3.13}}, is left invariant on $(\widehat M,\widehat
\Delta)$.
\end{Proposition}

\begin{proof}
Take $y\in \mathcal N_{\widehat\varphi}$.
Take also a~vector $\xi\in\mathcal H$ and let $\omega=\langle\,\cdot\,\xi,\xi\rangle$.
Consider an orthonormal basis $(\xi_i)$ in $\mathcal H$.
Then we have
\begin{gather*}
(\omega\otimes\iota)\widehat\Delta(y^*y)=\sum y_i^* y_i,
\end{gather*}
where $y_i=(\langle\,\cdot\,\xi,\xi_i\rangle \otimes\iota)\widehat\Delta(y)$.
We know from the previous proposition that $y_i\in\mathcal N_{\widehat\varphi}$ and that
$\Lambda_{\widehat\varphi}(y_i)=z_i\Lambda_{\widehat\varphi}(x)$ where $z_i=(\langle\,\cdot\,\xi,\xi_i\rangle
\otimes\iota)\widehat W^*$.
And because
\begin{gather*}
\sum z_i^*z_i=(\langle\,\cdot\,\xi,\xi\rangle\otimes\iota)(\widehat W\widehat W^*)=\omega(1)1,
\end{gather*}
it follows that
\begin{gather*}
\widehat\varphi((\omega\otimes\iota)\widehat \Delta(y^* y))  = \sum \widehat\varphi(y_i^* y_i)
 = \sum \langle z_i^*z_i \Lambda_{\widehat\varphi}(y),\Lambda_{\widehat\varphi}(y) \rangle
 = \omega(1)\widehat\varphi(y^*y).
\end{gather*}
This proves invariance.
\end{proof}

Observe that the invariance is proven by f\/irst constructing the candidate for the left regular representation and using
that this is a~unitary.
This is a~standard technique (e.g.\ see the construction of the Haar weight in~\cite{VD6}).

Now we are almost ready for the {\it main result}.
We just need that there is also a~right invariant Haar weight.
This will follow from the next proposition (a result that was already announced earlier).

\begin{Proposition}
Define $\widehat R$ on $\widehat M$ by $\widehat R(y)=Jy^* J$, where as before,~$J$ is the modular conjugation
associated with the weight~$\varphi$ on~$M$ $($cf.\ Notation~{\rm \ref{1.19})}.
Then $\widehat R$ is an involutive $*$-anti-automorphism of $\widehat M$ that flips the coproduct~$\widehat \Delta$.
\end{Proposition}

This all essentially follows from the formula i) in Proposition~\ref{1.20} (compare also with the proof of iii) in Theorem~\ref{2.3}).

As an immediate consequence, we get that the weight $\widehat\psi$, def\/ined as the composition of $\widehat \varphi$
with $\widehat R$, will be a~right invariant weight.
Therefore we have completed the proof of the following, main result.

\begin{Theorem}
The pair $(\widehat M,\widehat \Delta)$ $($as constructed in~Definitions~{\rm \ref{3.1}} and~{\rm \ref{3.5})}, is a~locally compact quantum group
$($in the sense of Definition~{\rm \ref{2.1})}.
\end{Theorem}

In the next section we will argue that the involutive $*$-anti-automorphism $\widehat R$ is indeed the unitary antipode
on the dual and that this notation is consistent.
In fact, we will obtain more formulas in the next section relating the objects of the original quantum group
$(M,\Delta)$ with those of the dual $(\widehat M,\widehat\Delta)$.

\medskip

{\bf The bidual.}
We f\/inish this section with some remark about biduality.
We have the following result.

\begin{Theorem}
The dual of $(\widehat M,\widehat \Delta)$ is again $(M,\Delta)$.
\end{Theorem}

This follows from the fact that the regular representation $\widehat W$ of the dual coincides with $\Sigma W^*\Sigma$ (cf.\ Proposition~\ref{3.14}).

Let us now make a~comparison with the theory of (multiplier) Hopf algebras.
There the dual is usually equipped with the coproduct, dual to the product and not with the opposite coproduct as we have done here.
It is obvious that in that case, the dual of the dual is the original algebra with the original coproduct.
If, as is done here in the operator algebra approach, the dual is equipped with the opposite coproduct, then one might
expect that the dual of the dual will yield the original algebra but with both the opposite product and the opposite coproduct.
That this is not seen here is simply a~result of the treatment.
There is no problem as the unitary antipode is a~map that converts the product to the opposite product and the coproduct
to the opposite coproduct.

\section{A collection of formulas}\label{Section5}

In this section, we collect many of the {\it formulas relating the various objects} associated with a~locally compact
quantum group and its dual.
We will not give all the possible relations (as there are many), but the most important ones.
Other equalities can easily be obtained from the ones that we prove.
Again it should be mentioned that we do not get really new results but that some of the results are proven in a~slightly
other fashion than in the original papers by Kustermans and Vaes.
Also the formulas are organized in another (perhaps more systematic) way.

Fix a~locally compact quantum group $(M,\Delta)$ and consider the dual $(\widehat M, \widehat \Delta)$ as constructed in
the previous section.
In this section we will freely use the def\/initions and notations of the previous sections.
When appropriate we will explicitly recall the necessary notions and results.

We will do so with the following def\/inition and notations (some of them introduced already in the proof of Theorem~\ref{2.10}).

\begin{Definition}\label{4.1}
Def\/ine continuous one-parameter groups of unitaries $(u_t)$, $(v_t)$ and $(w_t)$ on $\mathcal H$ by
\begin{gather*}
u_t\Lambda_\varphi(x)  = \Lambda_\varphi(\sigma_t(x)),
\qquad
v_t\Lambda_\varphi(x)  = \nu^{\frac12 t} \Lambda_\varphi(\tau_t(x)),
\qquad
w_t\Lambda_\varphi(x)  = \nu^{-\frac12 t} \Lambda_\varphi(\sigma'_t(x)),
\end{gather*}
when $x\in \mathcal N_\varphi$.
\end{Definition}

The relative invariance of~$\varphi$ with respect to $(\tau_t)$ and $(\sigma'_t)$ respectively is used to justify the
def\/initions of $(v_t)$ and $(w_t)$.
Recall from Theorems~\ref{2.7} and~\ref{2.9} that $\varphi\circ\tau_t=\nu^{-t} \varphi$ and
$\varphi\circ\sigma'_t =\nu^{t} \varphi$ (where~$\nu$ is the scaling constant).
Of course $u_t=\nabla^{it}$ for all~$t$.
We have just introduced this notation in order to have some more symmetry.
In what follows, we will use either of these two notations for this one-parameter group.

\begin{Remark}\quad
\begin{enumerate}\itemsep=0pt
\item[i)] We have that $v_t=P^{it}$ where~$P$ is the operator def\/ined in~\cite[Def\/inition~6.9]{K-V2}.
Because of the special role of this operator (see further), and because we want to be as close as possible to the
notations used in the papers by Kustermans and Vaes, we will further in this section use~$P^{it}$ as well as~$v_t$
(whatever is more convenient), just as in the case of~$\nabla^{it}$ and~$u_t$.

\item[ii)] It can be verif\/ied that $w_t\Lambda_\psi(x)=\Lambda_\psi(\sigma_t'(x))$ for all $x\in \mathcal N_\psi$.
Therefore, $w_t={\nablad}^{it}$ where~$\nablad$ is the modular operator associated with the right Haar weight~$\psi$
on~$M$.
Also in this case we will use ${\nablad}^{it}$ as well as $w_t$.
\end{enumerate}
\end{Remark}

In the following proposition we formulate a~f\/irst relation involving some of these operators.

\begin{Proposition}\label{4.3}
We have $\nablad^{it}=\delta^{it}
(J\delta^{it}J)
\nabla^{it}$ for all~$t$.
\end{Proposition}

This result follows because $\sigma'_t(x)=\delta^{it} \sigma_t(x) \delta^{-it}$ for all $x\in M$ and
$\sigma_s(\delta^{it})=\nu^{ist}\delta^{it}$ (see Theorem~\ref{2.11}).
Indeed, take $x\in M$ and assume that $x\delta^{\frac12}$ is bounded and belongs to $\mathcal N_\varphi$.
Then the following calculation is justif\/ied.
We have{\samepage
\begin{gather*}
\begin{split}
& {\nablad}^{it}\Lambda_\psi(x) =\Lambda_\psi(\sigma'(x))=\Lambda_\varphi\big(\sigma'_t(x)\delta^{\frac12}\big)
 =\nu^{-\frac12 t}\Lambda_\varphi\big(\sigma'_t\big(x\delta^{\frac12}\big)\big)
 =\nu^{-\frac12 t}\Lambda_\varphi\big(\delta^{it}\sigma_t\big(x\delta^{\frac12}\big)\delta^{-it}\big)
\\
& \phantom{{\nablad}^{it}\Lambda_\psi(x)}
 =\nu^{-\frac12 t}\delta^{it}J\sigma_{-\frac{i}{2}}\big(\delta^{it}\big)J\Lambda_\varphi\big(x\delta^{\frac12}\big)
 =\nu^{-\frac12 t}\nu^{\frac12 t}\delta^{it}J\delta^{it}J\Lambda_\varphi\big(x\delta^{\frac12}\big)
 \end{split}
\end{gather*}
and we get the result.}

Observe that $\delta^{it}$ and $J\delta^{it}J$ commute and that also $\nabla^{it}$ commutes with the product
$\delta^{it} (J\delta^{it}J)$.

From the def\/initions of these unitaries and using that the automorphism groups involved mutually commute (cf.\
Theorem~\ref{2.8}), we also get the following formulas.

\begin{Proposition}\label{4.4}\quad
\begin{enumerate}\itemsep=0pt
\item[$i)$] All the unitaries $(u_t)$, $(v_t)$ and $(w_t)$ mutually commute and they all also commute with the modular
conjugation~$J$.

\item[$ii)$] We have
\begin{gather*}
\sigma_t(x)  = u_t x u_t^*,
\qquad
\sigma'_t(x)  =w_t x w_t^*,
\qquad
\tau_t(x)  = v_t x v_t^*
\end{gather*}
for all $x\in M$ and $t\in {\mathbb R}$.
\end{enumerate}
\end{Proposition}

We will now introduce some more operators of the same type as in Def\/inition~\ref{4.1}.
Later we will justify the notations used in this def\/inition.

\begin{Definition}\label{4.5}
Def\/ine a~conjugate linear, involutive operator $\widehat J$ and a~non-singular, positive self-adjoint operator $\widehat
\nabla$ on $\mathcal H$~by
\begin{gather*}
\widehat J \Lambda_\varphi(x)  = \Lambda_\psi(R(x)^*),
\qquad
\widehat\nabla^{it} \Lambda_\varphi (x)  = \Lambda_\varphi \big(\tau_t(x)\delta^{-it}\big),
\end{gather*}
whenever $x\in \mathcal N_\varphi$ and $t\in {\mathbb R}$.
\end{Definition}

First recall the following from the remark made at the end of Section~\ref{Section3}.
If $x\in M$ and $R(x)^*\delta^{\frac12}$ is bounded and belongs to $\mathcal N_\varphi$, then $R(x)^*\in \mathcal
N_\psi$ and
\begin{gather*}
\Lambda_\psi(R(x)^*)=\Lambda_\varphi\big(R(x)^*\delta^{\frac12}\big).
\end{gather*}
In fact, for any $x\in M$ we have that $x\in \mathcal N_\varphi$ if and only if $R(x)^*\in \mathcal N_\psi$ because
$\psi (R(x^*x)) = \varphi(x^*x)$.
It follows that $\widehat J$ is well-def\/ined and that it is isometric.
Because $R(\delta)=\delta^{-1}$, we will have that $\widehat J^2=1$.

Also $\tau_t(x) \delta^{-it}\in \mathcal N_\varphi$ whenever $x\in \mathcal N_\varphi$.
To show that the map $\Lambda_\varphi(x) \mapsto \Lambda_\varphi(\tau_t(x)\delta^{-it})$ is isometric and that we get
indeed a~one-parameter group of unitaries, we can either make a~straightforward calculation or use the formula that we
obtain in Proposition~\ref{4.6} below.

It is also straightforward to show that $\widehat J$ and $\widehat\nabla^{it}$ commute for all $t\in {\mathbb R}$.

From the def\/inition, we immediately get the following relation.

\begin{Proposition}\label{4.6}
We have $\widehat\nabla^{it}=J\delta^{it}J
P^{it}$ for all~$t$.
\end{Proposition}

This formula is an easy consequence of the def\/initions and again of the fact that
$\sigma_s(\delta^{it})=\nu^{ist}\delta^{it}$.
Observe that also here the operators $J\delta^{it}J$ and $P^{it}$ commute.

\begin{Proposition}\label{4.7}
We have $R(x)=\widehat J x^* \widehat J$ and $\tau_t(x)=\widehat\nabla^{it} x \widehat\nabla^{-it}$ for all~$x$ and all~$t$.
\end{Proposition}

Again these two formulas are easy consequences of the def\/initions.
Remark that here we get another one-parameter group of unitaries that implements the scaling group $\tau_t$ (compare with Proposition~\ref{4.4}).

Recall that we also had the formulas $R(x)=Ix^*I$ and $\tau_t(x)=L^{it}x L^{-it}$ (see Def\/inition~\ref{1.22}).
Indeed, there are reasons to believe that we have $I=\widehat J$ and $L=\widehat\nabla$.
We will come back to this problem later (cf.\ Remark~\ref{4.10}).

Now we will prove some new relations.

\begin{Proposition}\label{4.8}
We have $\nablad^{it}=\widehat J \nabla^{-it} \widehat J$ for all~$t$.
\end{Proposition}

The result follows easily from the fact that $R(\sigma_t(x))=\sigma'_{-t}(R(x))$ for all $x\in M$ (cf.\
Theorem~\ref{2.6}) and $\sigma_t(\delta)=\nu^t\delta$ (Theorem~\ref{2.11}).

This is one useful formula involving the operator $\widehat J$.
Another one is the following relation between the left regular representation~$W$ associated with~$\varphi$ and the
right regular representation~$V$ associated with~$\psi$ on~$M$.

\begin{Proposition}
We have
$V = \big(\widehat J\otimes \widehat J\,\big)\Sigma W^* \Sigma\big(\widehat J \otimes \widehat J\,\big)$.
\end{Proposition}

Recall that here~$\Sigma$ denotes the f\/lip on $\mathcal H \otimes \mathcal H$.

The proof is straightforward.
Formally we can write, with the Sweedler notation $\Delta(x)=x_{(1)} \otimes x_{(2)}$ (without the summation sign
because we are using this symbol for something else here) and using that~$R$ f\/lips the coproduct:
\begin{gather*}
V\big(\widehat J \otimes \widehat J\,\big) (\Lambda_\varphi(x) \otimes \xi)  = V\big(\Lambda_\varphi\big(R(x)^*\delta^\frac12\big) \otimes \widehat J \xi\big)
 = V\big(\Lambda_\psi(R(x)^*) \otimes \widehat J \xi\big)
\\
\phantom{V\big(\widehat J \otimes \widehat J\,\big) (\Lambda_\varphi(x) \otimes \xi)}
 = \Lambda_\psi(R(x_{(2)})^*) \otimes R(x_{(1)})^* \widehat J \xi
 = \Lambda_\varphi\big(R(x_{(2)})^* \delta^\frac12\big) \otimes \widehat J x_{(1)} \xi
\\
\phantom{V\big(\widehat J \otimes \widehat J\,\big) (\Lambda_\varphi(x) \otimes \xi)}
 = \big(\widehat J \otimes \widehat J\,\big)(\Lambda_\varphi(x_{(2)}) \otimes x_{(1)} \xi),
\end{gather*}
when $x\in \mathcal N_\varphi$ and $\xi \in \mathcal H$.
We see that indeed $(\widehat J\otimes \widehat J)V(\widehat J \otimes \widehat J) = \Sigma W^* \Sigma$.

As we mentioned already, we will show later that $\widehat J$ is the modular conjugation associated with
$\widehat\varphi$ on $\widehat M$ and therefore we will get $\widehat J\widehat M \widehat J = \widehat M'$ and $V\in \widehat M' \otimes M$.

It would be possible to include more relations at this point, but we will postpone this.
First we will show that indeed $\widehat J$ and $\widehat \nabla$ are the modular conjugation and the modular operator of $\widehat\varphi$.

In Section~\ref{Section2}, we have mentioned that formally $K\Lambda_\psi(x)=\Lambda_\psi(S(x)^*)$ for well chosen elements $x\in\mathcal N_\psi$.
From the proof of Lemma~\ref{3.12} we also expect that $K^*\Lambda_\varphi(x)=\Lambda_\varphi(S(x^*))$ for certain
elements $x\in \mathcal N_\varphi$.
We have mentioned these formulas only for a~better understanding and motivation.
We did not use these formulas in any argument.

Now, as promised, we want to look at these formulas in a~correct way.
Remember that we use $K=IL^{\frac12}$ to denote the polar decomposition of~$K$ and that $K^*=IL^{-\frac12}$.
As we have seen in the previous section (see Lemma~\ref{3.11}), we expect~$K$ to be the closure of the map
$\Lambda_{\widehat\varphi}(y) \mapsto \Lambda_{\widehat\varphi}(y^*)$ for $y\in \mathcal N_{\widehat\varphi} \cap
\mathcal N_{\widehat\varphi}^*$.

We will prove a~closely related result.
It is not quite the same as we explain in the following important remark.

\begin{Remark}\label{4.10}\quad
\begin{enumerate}\itemsep=0pt
\item[i)] Formally we have
\begin{gather*}
\widehat J \widehat\nabla^{\frac12} \Lambda_\varphi(x) = \widehat J \Lambda_\varphi\big(\tau_{-\frac{i}{2}}(x)\delta^{-\frac12}\big)
 = \Lambda_\varphi \big(R\big(\tau_{-\frac{i}{2}}(x)\big)^* R\big(\delta^{-\frac12}\big)^* \delta^\frac12\big)
 = \Lambda_\varphi(S(x)^*\delta)
\end{gather*}
and
\begin{gather*}
\widehat J \widehat\nabla^{-\frac12}\Lambda_\varphi(x) = \widehat J \Lambda_\varphi\big(\tau_{\frac{i}{2}}(x)\delta^{\frac12}\big)
 = \Lambda_\varphi \big(R\big(\tau_{\frac{i}{2}}(x)\big)^* R\big(\delta^{\frac12}\big)^* \delta^\frac12\big)
 = \Lambda_\varphi(S(x^*))
\end{gather*}
for appropriate elements $x\in \mathcal N_\varphi$.
Therefore we expect that indeed $K=\widehat J\widehat\nabla^\frac12$ and so $I=\widehat J$ and $L=\widehat\nabla$.
However, although this result is likely to be true, it is not clear how to prove~it.

\item[ii)] We will not worry about this question (in this paper).
After all we only have used~$I$ and~$L^{it}$ in the formulas $R(x)=Ix^*I$ and $\tau_t(x)=L^{it}x L^{-it}$ and these
formulas remain true when~$I$ and~$L$ are replaced by $\widehat J$ and $\widehat\nabla$ (cf.\ Proposition~\ref{4.7}).
In other words, we can safely replace~$I$ by~$\widehat J$ and~$L$ by~$\widehat\nabla$ in the relevant formulas.

\item[iii)] We plan to look at this problem closer in the notes~\cite{VD8}.
See also~\cite{K-V3} but remark again that the operator~$K$ here should be compared with the operator $G^*$ (see e.g.\
the introduction and Corollary~2.9 in~\cite{K-V3}).
\end{enumerate}
\end{Remark}

Now we come to the following important result.

\begin{Theorem}\label{4.11}
Denote by $\widehat T$ the closure of the conjugate linear map $\Lambda_{\widehat\varphi}(y)\mapsto
\Lambda_{\widehat\varphi}(y^*)$ where $y \in \mathcal N_{\widehat\varphi} \cap \mathcal N_{\widehat\varphi}^*$.
Then the polar decomposition of $\widehat T$ is given by $\widehat T=\widehat J \widehat\nabla^\frac12$ $($with the
operators $\widehat J$ and $\widehat\nabla$ as in Definition~{\rm \ref{4.5})}.
\end{Theorem}
\begin{proof}
From Lemma~\ref{3.11} we know that
\begin{gather*}
\langle\Lambda_{\widehat\varphi}(y^*), \Lambda_{\varphi}(a) \rangle
=\langle\Lambda_{\widehat\varphi}(y),\Lambda_{\varphi}(S(a^*))\rangle^-,
\end{gather*}
where $y\in \widehat{\mathcal N} \cap \widehat{\mathcal N}^*$ and when $a\in \mathcal N_\varphi$ is an element such that
also $a^*\in\mathcal D(S)$ and $S(a^*)$ is still in $\mathcal N_\varphi$.
It follows that $\widehat T$ is contained in the adjoint of the map $\Lambda_\varphi(a) \mapsto \Lambda_\varphi(S(a^*))$ with~$a$ as above.

As we have seen already in the proof of Lemma~\ref{3.12}, we can relatively easily produce such elements $a\in \mathcal
N_\varphi$ by requiring that~$a$ is analytic both with respect to the automorphism group~$\tau$, as well as with respect
to multiplication from the right with the unitary group $(\delta^{it})$.
Moreover the space of elements $\Lambda_\varphi(a)$ with such elements~$a$ will be left invariant by the operators
$\widehat\nabla^{it}$.
All of this will imply that $\Lambda_{\widehat\varphi}(y) \in \mathcal D(\widehat J \widehat\nabla^\frac12)$ and that
$\widehat J\widehat\nabla^\frac12 \Lambda_{\widehat\varphi}(y)= \Lambda_{\widehat\varphi}(y^*)$ whenever $y \in
\widehat{\mathcal N} \cap \widehat{\mathcal N}^*$.
Now, a~straightforward calculation shows that $\Lambda_{\widehat\varphi}(\widehat{\mathcal N} \cap \widehat{\mathcal
N}^*)$ is invariant under the unitaries $\widehat\nabla^{it}$.
Then the result follows.
\end{proof}

So, we get as expected, that $\widehat J$ {\it is the modular conjugation} and $\widehat\nabla$ the {\it modular
operator} associated with the dual left Haar weight $\widehat\varphi$.

As a~f\/irst consequence of this result, we get e.g.\ that $M\cap \widehat M={\mathbb C} 1$.
Indeed, just before Lemma~\ref{3.4} in the previous section, we saw that $M\cap \widehat M'={\mathbb C} 1$.
And because $\widehat J M \widehat J=M$ and $\widehat J \widehat M \widehat J = \widehat M'$, we get that also $M\cap
\widehat M={\mathbb C} 1$.

A~second
important consequence is the following result (see~\cite[Proposition~2.1]{K-V3}).

\begin{Proposition}\label{4.12}
The unitary antipode $\widehat R$ on $\widehat M$ is given by $\widehat R(y)=Jy^*J$.
The scaling group~$\widehat \tau$ is given by~$\widehat \tau_t(y)=\nabla^{it} y \nabla^{-it}$ for all $t\in {\mathbb R}$
whenever $y\in\widehat M$.
\end{Proposition}

This result can be proven in two ways.
One argument uses duality as follows.
We know that the unitary antipode~$R$ on~$M$ is given by $R(x)=\widehat Jx^*\widehat J$ and the scaling group by
$\tau_t(x)=\widehat\nabla^{it} x \widehat\nabla^{-it}$ (cf.\ Proposition~\ref{4.7}).
From Theorem~\ref{4.11}, we know that $\widehat T = \widehat J \widehat\nabla^\frac12$ is the polar decomposition of the
dual operator $\widehat T$.
Hence, because $J \nabla^\frac12$ is the polar decomposition of the operator~$T$ associated with the left Haar
weight~$\varphi$ on~$M$, we will get the formulas in the proposition by duality.
A~second argument would be possible by using the formulas
\begin{gather*}
\big(\widehat J \otimes J\big) W \big(\widehat J \otimes J\big) = W^*
\qquad
\text{and}
\qquad
\big(\widehat\nabla^{it} \otimes \nabla^{it}\big) W \big(\widehat\nabla^{-it} \otimes \nabla^{-it}\big) = W
\end{gather*}
and the fact that the right leg of~$W$ is dense in $\widehat M$ (and related things).
Remark that these formulas are found in Proposition~\ref{1.20} with the operators~$I$ and~$L$ in the place of $\widehat J$
and~$\widehat\nabla$ but, as we explained in Remark~\ref{4.10}ii), they will still be correct.

Next we have another important consequence.

\begin{Proposition}\label{4.13}
The scaling constant $\widehat \nu$ of the dual is $\nu^{-1}$.
Furthermore, we also have
\begin{gather*}
P^{it} \Lambda_{\widehat\varphi}(y) = \nu^{-\frac12 t} \Lambda_{\widehat\varphi}(\widehat\tau_t(y))
\end{gather*}
for all $y\in \widehat M$.
\end{Proposition}
\begin{proof}
The proof of this result will use the basic formula $(\tau_t\otimes\widehat\tau_t)W = W$ for all~$t$.

Start with an element $y\in \widehat{\mathcal N}$ with $y=(\omega\otimes\iota) W$ and $\omega\in M_*$ so that
$\omega(a^*) = \langle \widehat\Lambda(y), \Lambda_\varphi(a)\rangle$ for all $a\in \mathcal N_\varphi$.
Then
\begin{gather*}
\langle P^{it} \widehat\Lambda(y), \Lambda_\varphi(a) \rangle = \langle \widehat\Lambda(y), P^{-it}\Lambda_\varphi(a) \rangle
 = \nu^{-\frac12 t}\langle \widehat\Lambda(y), \Lambda_\varphi(\tau_{-t}(a)) \rangle
 = \nu^{-\frac12 t} \omega(\tau_{-t} (a^*))
\end{gather*}
for all $a\in \mathcal N_\varphi$.
Now we have $((\omega\circ\tau_{-t}) \otimes \iota) W = \widehat\tau_t(y)$ by the basic formula above.
Hence, we see that $\widehat\tau_t(y)\in\widehat{\mathcal N}$ and that $P^{it} \widehat\Lambda(y) = \nu^{-\frac12 t}
\widehat\Lambda(\widehat\tau_t(y))$.
Then this formula (as in the formulation of the proposition) holds for all $y\in \mathcal N_{\widehat\varphi}$.

From this formula (and because $P^{it}$ is a~unitary operator), it follows that $\widehat\varphi \circ \widehat\tau_t =
\nu^t \widehat\varphi$ for all~$t$ and therefore $\widehat\nu = \nu^{-1}$.
\end{proof}

We see that $(P^{it})$ is, in some sense, a~self-dual group of unitaries.
One can write $\widehat P= P$ when $\widehat P$ would be def\/ined for the dual as~$P$ is def\/ined for the original pair
$(M,\Delta)$.
So, this one-parameter group of unitaries implements the scaling group both on~$M$ and on $\widehat M$.
In combination with the formula $(S \otimes \iota)W=W^*$, it will give that~$W$ is a~manageable multiplicative unitary
(in the sense of~\cite{W2}).
As we have already mentioned however, we will not use this property as such.

We will now use a~similar technique as in the proof of Proposition~\ref{4.13} to obtain more relations.
Recall that in the above proof, we used the basic formula $(\tau_t\otimes\widehat\tau_t)W = W$.
First, we will now formulate again some other results of this type.

Recall that formally we have $W^*(\xi\otimes\Lambda_\varphi(x))=\sum x_{(1)}\xi \otimes \Lambda_\varphi(x_{(2)})$ (when
we use $\Delta(x)=\sum x_{(1)}\otimes x_{(2)}$) with $x\in \mathcal N_\varphi$ and $\xi \in \mathcal H$.
Then from the def\/initions of $(u_t)$, $(v_t)$ and $(w_t)$ and the formulas
\begin{alignat*}{3}
& \Delta(\sigma_t(x))  = (\tau_t \otimes \sigma_t) \Delta(x),
\qquad &&
\Delta(\sigma'_t(x))  = (\sigma'_t \otimes \tau_{-t}) \Delta(x),&
\\
& \Delta(\tau_t(x))  = (\tau_t \otimes \tau_t) \Delta(x),
\qquad &&
\Delta(\tau_t(x))  = (\sigma_t \otimes \sigma'_{-t}) \Delta(x),&
\end{alignat*}
for all $x\in M$ and $t\in {\mathbb R}$, we get the following relations (see also the proof of Theorem~\ref{2.10}).

\begin{Lemma}\label{4.14}
We have
\begin{alignat*}{3}
& (\tau_t \otimes \iota) W  = (1\otimes u_t^*)W(1\otimes u_t),
\qquad &&
(\sigma'_t \otimes \iota) W  = (1\otimes w_t^*)W(1\otimes v_t^*), &
\\
& (\tau_t \otimes \iota) W  = (1\otimes v_t^*)W(1\otimes v_t),
\qquad &&
(\sigma_t \otimes \iota) W  = (1\otimes v_t^*)W(1\otimes w_t^*)&
\end{alignat*}
for all $t\in {\mathbb R}$.
\end{Lemma}

Each of the four formulas comes from one of the above relations, in the same order.
The f\/irst and the third formula do not give anything new.
This is simply the fact that $(\tau_t\otimes\widehat\tau_t)W=W$, combined with the knowledge that both $u_t$ and $v_t$
implement $\widehat\tau_t$ on $\widehat M$.

However, using all these results, in combination with the def\/inition of the dual map $\widehat \Lambda$, as in the proof
of Proposition~\ref{4.13}, we obtain the following result about the modular element $\widehat\delta$, relating the left
and the right Haar weights on the dual $(\widehat M,\widehat\Delta)$.

\begin{Proposition}\label{4.15}
We have $\widehat \delta^{it} = v_t^*w_t^* $ where $v_t$ and $w_t$ are defined as in Definition~{\rm \ref{4.1}}.
Furthermore $\widehat\Delta (\widehat\delta)^{it} = \widehat\delta^{it} \otimes \widehat\delta^{it}$ for all~$t$.
\end{Proposition}
\begin{proof}
Because we have
\begin{gather*}
(\sigma_t \otimes\iota) W  = (1 \otimes v_t^*) W (1 \otimes w_t^*)
 = (1 \otimes v_t^*) W (1 \otimes v_t)(1 \otimes v_t^* w_t^*)
 = ((\tau_t\otimes\iota) W) (1 \otimes v_t^* w_t^*),
\end{gather*}
we see that $w_tv_t\in\widehat M$ for all~$t$.

Now if we use the same type of argument as in the proof of the previous proposition, we get $\widehat\tau_{t}(y) v_{t}
w_{t} \in \widehat{\mathcal N}$ when $y\in \widehat{\mathcal N}$ and
\begin{gather*}
\nabla^{it} \widehat\Lambda(y) = \widehat\Lambda(\widehat\tau_{t}(y) v_{t}w_{t}).
\end{gather*}
If we compare this with the earlier and similar formula
\begin{gather*}
\widehat\nabla^{it} \Lambda_\varphi(y) = \Lambda_\varphi\big(\tau_{t}(x) \delta^{-it}\big)
\end{gather*}
for $x\in \mathcal N_\varphi$, we see that we must have $\widehat\delta^{it}=v_t^* w_t^*$ for all~$t$ (by duality).
Recall that $v_t$ and $w_t$ commute with each other.

From the formula $(\sigma'_t \otimes \iota) W = (1\otimes w_t^*) W (1 \otimes v_t^*)$ we get
\begin{gather*}
(w_t \otimes w_t) W (w_t^* \otimes v_t) = W
\end{gather*}
and from the formula $(\tau_t \otimes \iota) W = (1\otimes v_t^*) W (1 \otimes v_t)$ we f\/ind
\begin{gather*}
(v_t \otimes v_t) W (v_t^* \otimes v_t^*) = W
\end{gather*}
and combining these two results, we obtain
\begin{gather*}
\big(\widehat\delta^{-it} \otimes \widehat\delta^{-it}\big) W \big(\widehat\delta^{it} \otimes 1\big) = W
\end{gather*}
proving that $\widehat\Delta(\widehat\delta^{it}) = \widehat\delta^{it} \otimes \widehat\delta^{it}$.
\end{proof}

\begin{Remark}\label{4.16}
By duality, we also will get $\Delta(\delta^{it})=\delta^{it} \otimes \delta^{it}$ for all~$t$.

This seems to be a~strange (and certainly not an obvious) way to prove this basic formula.
One could expect a~more direct proof (e.g.\ along the lines of the proof of this formula for algebraic quantum groups in~\cite{VD4}).
However, an attempt to do this turns out to become quite involved (see  in~\cite[Section~7]{K-V2}).
In~\cite{M-N-W}, this formula is proven in a~more elegant way, but also uses results, not only about $(M,\Delta)$
itself, but also about the dual $(\widehat M,\widehat\Delta)$ (see  in~\cite[Proposition~6.12]{M-N-W}).
\end{Remark}

Now we are ready to collect the main formulas and relations.
First, we have a~number of formulas that express some of the operators in terms of the others.

\begin{Theorem}\label{4.17}
We have the following formulas for the four modular operators:
\begin{alignat*}{3}
& \nabla^{it} = \big(\widehat J \widehat \delta^{it} \widehat J\,\big)  P^{it},
\qquad && \nablad^{it} =\widehat\delta^{-it}P^{-it} =\widehat J \nabla^{-it} \widehat J,&\\
&
\widehat\nabla^{it} = \big(J \delta^{it} J\big)  P^{it},
\qquad && \widehat\nablad^{it} = \delta^{-it}P^{-it}=J \widehat\nabla^{-it} J.&
\end{alignat*}
\end{Theorem}

The second formula was obtained in Proposition~\ref{4.6}.
The f\/irst one follows e.g.\ by duality.
From Proposition~\ref{4.8} we get the third formula while the last one follows again by duality.

We also have a~number of {\it commutation rules}.
We know that $P^{it}$ commutes with~$J$, $\widehat J$,~$\delta$ and~$\widehat\delta$ (and as a~consequence of the above
relations, also with all these modular operators).
We also know that~$J$ commutes with $\widehat \delta^{it}$ (because $\widehat R(\widehat\delta)=\widehat\delta^{-1}$)
and that $\widehat J$ commutes with $\delta^{it}$ (because $R(\delta)=\delta^{-1}$).
The non-trivial commutation rules are formulated in the next two theorems.

\begin{Theorem}\label{4.18}
We have
\begin{gather*}
\widehat J \nabla^{-it} \widehat J  = \delta^{it} (J\delta^{it}J) \nabla^{it},
\qquad
J \widehat\nabla^{-it} J  = \widehat\delta^{it}\big(\widehat J\widehat\delta^{it}\widehat J\,\big) \widehat\nabla^{it}.
\end{gather*}
\end{Theorem}

The f\/irst formula is a~combination of the formulas in Propositions~\ref{4.3} and~\ref{4.8}.
The second one comes with duality.

Next, we get the following two basic commutation rules.

\begin{Theorem}\label{4.19}
We have $\widehat\delta^{is}\delta^{it}=\nu^{-ist}\delta^{it}\widehat\delta^{is}$ for all $s,t\in {\mathbb R}$.
Also $\widehat J J= \nu^{\frac{i}{4}} J \widehat J$.
\end{Theorem}
\begin{proof}
There are several ways to prove these formulas.

To prove the f\/irst one, we can start with $\widehat \delta^{it} = v_t^*w_t^*$ as obtained in Proposition~\ref{4.15}.
We know that $w_s$ implements $\sigma'_s$ and that $\sigma'_s(\delta^{it})=\nu^{ist}\delta^{it}$ (see
Theorem~\ref{2.11}).
We also know that $v_s$ implements $\tau_s$ and that $\tau_s(\delta^{it})=\delta^{it}$ (again see Theorem~\ref{2.11}).
Combining these results will yield the f\/irst formula of the theorem.

To prove the second formula one calculates (carefully) that
\begin{gather*}
\widehat J T \widehat J=\nu^{\frac{i}{4}} T \delta^{\frac12}J\delta^{-\frac12}J,
\end{gather*}
where~$T$ is the closure of the map $\Lambda_\varphi(x)\mapsto \Lambda_\varphi(x^*)$ where $x\in \mathcal N_\varphi \cap
\mathcal N_\varphi$ as before (see Notation~\ref{1.19}).
Then the uniqueness of the polar decompostion of~$T$ gives the second formula.
\end{proof}

The uniqueness of the polar decomposition of~$T$, as used above, will also give the f\/irst formula of Theorem~\ref{4.18}.
On the other hand, it would also be possible to use this formula to give a~more direct proof of the commutation rule
between $\widehat J$ and~$J$.

If we combine the formulas of Theorem~\ref{4.17} with the commutation rules in Theorems~\ref{4.18} and~\ref{4.19}, we get most (if not all) of the other commutation rules.
There is no need to include them all here (see e.g.\ Proposition~2.13 in~\cite{K-V3}).
Some of these formulas, e.g.\ the second one in Theorem~\ref{4.19} above, giving the commutation rules between~$J$ and~$\widehat J$ is nicely illustrated in~\cite{M-VD2}.

Let us now {\it finish this section} with another formula, expressing the operators $P^{it}$ in terms of the other operators.
It is an analytic version of Radford's formula for the 4th power of the antipode, proven in the case of
a~f\/inite-dimensional Hopf algebra in~\cite{R} and extended to the case of algebraic quantum groups in~\cite{D-VD}.

\begin{Theorem}
We have
\begin{gather*}
P^{-2it} = \delta^{it}  (J\delta^{it}J)  \widehat\delta^{it}  \big(\widehat J\widehat\delta^{it} \widehat J\,\big)
\end{gather*}
for all~$t$.
\end{Theorem}
\begin{proof}
From Theorem~\ref{4.18} we get
\begin{gather*}
\widehat J \nabla^{-it} \widehat J = \delta^{it} (J\delta^{it}J) \nabla^{it}
\end{gather*}
and if we replace in this formula $\nabla^{it}$ by $(\widehat J \widehat \delta^{it} \widehat J) P^{it}$ two times, we
arrive at the desired result.
\end{proof}

This is indeed one of the possible analytical versions of the formula, valid in the theory of f\/inite-dimensional Hopf
algebras
\begin{gather*}
S^4(a)= \delta^{-1} \big(\widehat \delta \triangleright a\triangleleft \widehat\delta^{-1}\big) \delta,
\end{gather*}
where $\triangleright$ and $\triangleleft$ are used to denote the canonical left and right actions of the dual algebra.
Remark that $P^{2it} \Lambda_\varphi(a) = \nu^t\Lambda_\varphi (\tau_{2t}(a))$ and that $\tau_{-i}= S^2$.
Also $\delta^{it}\,J \delta^{it} J \Lambda_\varphi(a) = \nu^{\frac12 t}\Lambda_\varphi(\delta^{it} a~\delta^{-it})$.
To explain the last part of the formula, one should observe that the operator $\widehat\delta^{it}$ is a~convolution
operator on~$M$, but on the Hilbert space level and similarly for $\widehat J \widehat\delta^{it} \widehat J$.
The f\/irst one is `left' convolution and the second one is `right' convolution.

We have now the essential formulas, all formulated in these four last theorems.
Also the commutation rules with the left regular representation, as formulated in Lemma~\ref{4.14}, are essentially taken care of.

One can also draw certain interesting conclusions from these formulas.
If e.g.\ the left and the right Haar weight on~$M$ coincide, that is when $\delta=1$, it will follow from the
commutation rules in Theorem~\ref{4.19} that the scaling constant~$\nu$ has to be one.
If both~$\delta$ and $\widehat\delta$ are trivial, it will follow from the last theorem that also $P=1$ and that the
scaling automorphisms are trivial.
From Theorem~\ref{4.17}, it will follow that all the modular automorphisms are trivial in this case and so that all Haar
weights must be traces.
Surely, there are other arguments for these statements, but it is nice to see how they follow from all these formulas.

\section{Conclusion and possible further research}\label{Section6}

The original motivation for developing the theory of locally compact quantum groups was the desire to generalize
Pontryagin's duality for locally compact abelian groups to the non-abelian case.
After various interesting, but partial solutions, a~f\/irst more or less satisfactory theory was developed in the early
seventies, independently by Kac and Vainermann~\cite{V-K} on the one hand and Enock and Schwarz (see~\cite{E-S}) on the other hand.
They developed the theory of Kac algebras with a~duality between objects of the same type, generalizing the duality
between locally compact abelian groups as studied by Pontryagin.
Kac algebras have von Neumann algebras as their underlying operator algebras.

However, later it turned out that some new examples (like the ones developed by Drinfel'd in the theory of quantum
groups and the quantum $SU_2(q)$-group discovered by Woronowicz in the eighties) did not f\/it into the theory of Kac
algebras.
For Kac algebras, the square of the antipode is the identity map while this is not the case for these newer examples.
It became clear that a~new, still more general theory had to be found.

At that the same time,
there was the widely accepted philosophy that the expected theory should be formulated in
a~$C^*$-algebraic framework.
The reason is obvious.
Locally compact groups are topological objects and so locally compact quantum groups should be developed in
a~$C^*$-algebraic setting.
This eventually led to the results by Kustermans and Vaes as found in~\cite{K-V1} and~\cite{K-V2}.
A~locally compact quantum groups in these papers is indeed a~pair of a~$C^*$-algebra with a~coproduct and such that Haar weights exist.

Even though the theory was developed in a~$C^*$-setting, it was still obvious to consider also the von Neumann algebraic
version as this was from the very beginning expected to be equivalent, as for the earlier theories.
Again, this work was done by Kustermans and Vaes in~\cite{K-V3}.

In principle, the von Neumann algebraic approach in~\cite{K-V3} is independent of the earlier $C^*$-algebraic one.
However, it is very hard to read the paper~\cite{K-V3} without~\cite{K-V2} at hand as the two are still highly interconnected.

In this paper, we have studied locally compact quantum groups from the very beginning in a~von Neumann algebraic setting.
We did not rely on the earlier obtained $C^*$-version.
The argument we gave in the appendix, showing how to construct a~von Neumann algebraic locally compact quantum group
form a~$C^*$-algebraic one, without the need to develop the theory f\/irst, is crucial for this approach.
In this way, we can obtain an overall treatment of the theory of locally compact quantum groups within the easier and
more tractable von Neumann algebra context.

However, we have to admit that also the present development still requires a~good knowledge and quite some experience
with the technicalities that are typically encountered in the Tomita--Takesaki
theory of left Hilbert algebras and its
relation with faithful normal semi-f\/inite weights on von Neumann algebras.

We f\/inish here with the following observation.
In the theory of locally compact quantum groups, the existence of the Haar weights is part of the assumptions.
At the moment of this writing, there is not yet a~theory with reasonable assumptions from which the existence of the
Haar weights would follow.
Only in the case of compact quantum groups (see e.g.~\cite{W1} and~\cite{M-VD1}) and for discrete quantum groups
(see~\cite{VD3}), it is possible to formulate assumptions that allow to prove the existence of the Haar weights.

On the other hand, as it turns out, in the examples, it is not so dif\/f\/icult to f\/ind the Haar weights. See e.g.~\cite{VD6}.
And in fact, there is more.
As can be seen from the discussions in~\cite{VD6}, see e.g.\ Section~\ref{Section6}, in a~certain sense, it is possible to write down
a~formula for the weight that has to be the left Haar weight whenever it makes sense, that is essentially, whenever the
formula def\/ines a~semi-f\/inite weight.
This might open a~path to a~possible theory of locally compact quantum groups where the existence of the Haar weights
follow from the assumptions.

\appendix

\section{Appendix. Other approaches to locally compact\\ quantum groups}\label{appendixA}

In this appendix, we will relate the von Neumann algebra approach (in particular, as it is treated in this paper) with
the $C^*$-algebra approach.
We will consider the papers by Kustermans and Vaes~\cite{K-V1} and~\cite{K-V2}, as well as the paper by Masuda, Nakagami
and Woronowicz~\cite{M-N-W}.
We will also brief\/ly make a~comparison with the earlier paper by Masuda and Nakagami~\cite{M-N}.

Such a~comparison is certainly interesting.
In fact, the relation with (the equivalent) $C^*$-algebra approach is not only interesting, it is also an important
feature of the theory.
Nevertheless, because of the scope of this paper, we will be very brief here.
In~\cite{VD7}, some more details are already found, but we refer to forthcoming papers for all the details about the
material in this appendix (e.g.\ for results about weights on $C^*$-algebras and for the theory of locally compact quantum groups, see~\cite{VD8}.

Now, we start with a~$C^*$-algebra~$A$ and a~comultiplication on~$A$.
Recall the def\/inition:

\begin{Definition}\label{A.1}
A~{\it comultiplication on a~$C^*$-algebra~$A$} is a~non-degenerate $*$-homomor\-phism~$\Delta$ from~$A$ to the
multiplier algebra $M(A\otimes A)$ of the minimal $C^*$-tensor product $A\otimes A$ of~$A$ with itself, satisfying
coassociativity $(\Delta\otimes\iota)\Delta=(\iota\otimes\Delta)\Delta$.
It is also assumed that slices $(\omega\otimes\iota)\Delta(a)$ and $(\iota\otimes\omega)\Delta(a)$ belong to~$A$ for all
$a\in A$ and $\omega\in A^*$.
\end{Definition}

Non-degenerate here means that $\Delta(A)(A\otimes A)$ is dense in $A\otimes A$.
Because of this condition, the $*$-homomorphisms $\Delta\otimes\iota$ and $\iota\otimes\Delta$ have unique extensions
(still denoted in the same way) to unital $*$-homomorphisms from $M(A\otimes A)$ to $M(A\otimes A\otimes A)$.
Therefore, coassociativity, as formulated above, has a~meaning.
Slice maps are def\/ined from $M(A\otimes A)$ to $M(A)$.
So, in general these slices $(\omega\otimes\iota)\Delta(a)$ and $(\iota\otimes\omega)\Delta(a)$ belong to $M(A)$.
In~\cite{M-N-W}, it is assumed that $\Delta(A)(1\otimes A)$ and $\Delta(A)(A\otimes 1)$ are subsets of $A\otimes A$.
Because any $\omega\in A^*$ is of the form $\rho(a\,\cdot\,)$ (and of the form $\rho(\,\cdot\,a)$) with $\rho\in A^*$
and $a\in A$, the latter conditions will imply the conditions in the def\/inition above.
In fact, it is shown in the theory that also these stronger conditions are valid.
See e.g.\ the~remark following Proposition~\ref{1.21}.

The following result is not so dif\/f\/icult to obtain.

\begin{Proposition}\label{A.2}
Let $(M,\Delta)$ be a~locally compact quantum group $($as in Definition~{\rm \ref{2.1})} with left Haar
weight~$\varphi$.
Let~$W$ be the left regular representation $($as introduced in Section~{\rm \ref{Section2})} and define~$A$ to be the norm closure of the
space of slices $\{(\iota\otimes\omega)W \,|\, \omega\in \mathcal B(\mathcal H_\varphi)_* \}$.
Then~$A$ is a~$C^*$-algebra, it is a~$\sigma$-weakly dense subalgebra of~$M$ and the restriction of~$\Delta$ to~$A$ is
a~comultiplication on~$A$.
\end{Proposition}

It is easy to see (and a~standard result about multiplicative unitaries) that~$A$ is a~subalgebra.
And just as in the case of the dual (see Proposition~\ref{3.3}), one can show that it is
a~$*$-subalgebra.
As the left leg of~$W$ is dense in~$M$, it follows that also~$A$ is dense in~$M$.
To show that~$\Delta$ maps~$A$ into $M(A\otimes A)$ is again a~standard result for multiplicative unitaries
(see~\cite{B-S}).
See also the work on manageable multiplicative unitaries (see \cite{W2} and~\cite{S-W}).
Also here, we refer to~\cite{VD8} for a~detailed and independent approach.

Later in this appendix
(see Def\/inition~\ref{A.10}), we will recall the def\/inition of a~locally compact quantum group in
the $C^*$-algebra context, as given by Kustermans and Vaes in~\cite{K-V1} and~\cite{K-V2}.
First, we look at some results about (invariant) weights on $C^*$-algebras (with a~comultiplication).
The f\/irst one is about the kind of weights that is used in this theory.

\begin{Proposition}\label{A.3}
Let~$A$ be a~$C^*$-algebra.
Let~$\varphi$ be a~faithful, densely defined, lower semi-continuous weight on~$A$.
Denote by $\widetilde\varphi$ the normal weight on the double dual $A^{**}$ that extends~$\varphi$.
Denote by~$e$ the support projection of $\widetilde\varphi$ in $A^{**}$.
Then~$\varphi$ is approximately K.M.S.
$($see Definition~{\rm 1.34} in~{\rm \cite{K-V2})}
if and only if~$e$ is a~central projection in~$A^{**}$.
\end{Proposition}

Assume that~$\varphi$ is any faithful, densely def\/ined, lower semi-continuous weight on the $C^*$-algebra~$A$.
Consider its GNS representation $\pi_\varphi$ on the Hilbert space $\mathcal H_\varphi$.
Let $\Lambda_\varphi$ denote the canonical imbedding of $\mathcal N_\varphi$ in $\mathcal H_\varphi$.
A~vector $\xi$ in $\mathcal H_\varphi$ is called {\it right bounded} if $\Lambda_\varphi(a)\mapsto \pi_\varphi(a)\xi$,
where $a\in \mathcal N_\varphi$, is bounded.
In some sense, there are always plenty of right bounded vectors, but it can happen that the space of right bounded
vectors is not dense in $\mathcal H_\varphi$.
If we denote by~$q$ the projection onto the closure of the space of right bounded vectors, then~$q\in \pi_\varphi(A)''$.
It turns out to be the image of the support projection~$p$ of the unique normal extension of~$\widetilde\varphi$ to~$A^{**}$ under the unique normal extension $\widetilde {\pi_\varphi}$ of~$\pi_\varphi$ to~$A^{**}$.
Moreover, the support projection~$e$ of this extension $\widetilde {\pi_\varphi}$ is the central support of~$p$.

Recall that a~faithful, densely def\/ined lower semi-continuous weight on a~$C^*$-algebra is approximately K.M.S.
if it remains faithful when extended to $\pi_\varphi(A)''$.
This last property essentially means that the support of $\widetilde{\pi_\varphi}$ in $A^{**}$ (which is a~central
projection), coincides with the support projection of $\widetilde\varphi$.
This is the case if and only if the right bounded vectors are dense.
So, the above result should not come as a~surprise.
We refer to~\cite{VD8} for details.
Further in this appendix, we will speak about a~weight on a~$C^*$-algebra with central support, or shortly call it {\it a~central weight}.

Invariant weights on $C^*$-algebras with a~comultiplication are def\/ined as usual:

\begin{Definition}
Let~$A$ be a~$C^*$-algebra with a~comultiplication~$\Delta$ (as in Def\/inition~\ref{A.1}).
A~weight~$\varphi$ is called {\it left invariant} if $\varphi((\omega\otimes\iota)\Delta(a))=\|\omega\|\varphi(a)$
whenever $a\in A$, $a\geq 0$ and $\varphi(a)<\infty$ and when $\omega\in A^*$ and $\omega\geq 0$.
Similarly, a~{\it right invariant} weight is def\/ined.
\end{Definition}

Again, there is the following result.

\begin{Proposition}\label{A.5}
Let $(M,\Delta)$ be a~locally compact quantum group.
The restriction of the left Haar weight~$\varphi$ to the $C^*$-subalgebra~$A$ of~$M$, defined as in
Proposition~{\rm \ref{A.2}}, is a~faithful, densely defined, lower semi-continuous central and left invariant weight.
Similarly for the right Haar weight.
\end{Proposition}

\begin{proof}
It is quite obvious that these restrictions are faithful, lower semi-continuous and invariant weights.
They are central because they are restrictions of faithful weights to a~$C^*$-subalgebra of the von Neumann algebra~$M$.
To show that the right invariant weight~$\psi$ is still densely def\/ined on the $C^*$-algebra~$A$, one can use the result
in Lemma~\ref{1.13}.
To show that also the left invariant weight~$\varphi$ is still densely def\/ined on the $C^*$-algebra, one can use that
the unitary antipode~$R$ leaves the $C^*$-algebra invariant (a result that follows from the formula $(\widehat J\otimes
J)W(\widehat J\otimes J)=W^*$ and the def\/inition of~$R$).
\end{proof}

Using the formulas in the proof of Theorem~\ref{2.10} (or equivalently, the ones in Lemma~\ref{4.14}), we see that the
$C^*$-algebra~$A$ is invariant under the modular automorphism groups~$\sigma$ and $\sigma'$.
Therefore, the restrictions of the Haar weights~$\varphi$ and~$\psi$ are K-M-S weights (and so certainly central).
Observe also that the $C^*$-algebra is not only invariant under the unitary antipode, but also under the scaling
automorphisms~$\tau$, again see Lemma~\ref{4.14} (or equivalent earlier results).

We will now consider two important results.
They will enable us to go quickly from the $C^*$-algebra setting to the von Neumann algebra framework.

\begin{Proposition}\label{A.6}
Let~$A$ be a~$C^*$-algebra with a~comultiplication~$\Delta$.
Assume that~$\varphi$ is a~densely defined lower semi-continuous central weight on~$A$.
Consider the extension $\widetilde\varphi$ on $A^{**}$ of~$\varphi$.
Also extend~$\Delta$ to a~normal and unital $*$-homomorphism $\widetilde\Delta: A^{**} \to A^{**} \otimes A^{**}$ $($the
von Neumann algebra tensor product$)$.
Then $\widetilde\varphi$ is still left invariant.
\end{Proposition}

One uses left Hilbert algebra theory to show that the G.N.S.-space of the extension is the same as the original one.
This result is not completely trivial but is essentially proven in Section~2 of Chapter~VII on weights in~\cite{T3}.
Then it is shown that the extension is still invariant.
The way this is done, is similar as the invariance of the dual weight is proven in Proposition~\ref{3.15}.
Essentially, the argument used to prove that left invariance of~$\varphi$ implies the formula $WW^*=1$ for the left
regular representation~$W$, is `reversed'.

\begin{Proposition}\label{A.7}
As in Proposition~{\rm \ref{A.6}}, assume that~$A$ is a~$C^*$-algebra with a~comultiplication~$\Delta$.
Now assume that~$\varphi$ and~$\psi$ are non-trivial, densely defined, lower semi-continuous and central weights on~$A$
such that~$\varphi$ is left invariant and~$\psi$ is right invariant.
Then the supports of~$\varphi$ and~$\psi$ in $A^{**}$ are equal.
\end{Proposition}

\begin{proof}
Denote by~$e$ and~$f$ the supports of $\widetilde\varphi$ and $\widetilde\psi$ resp.
By assumption, they are central projections in $A^{**}$.
By the left invariance of $\widetilde\varphi$, we get $\widetilde\varphi((\omega\otimes\iota)\widetilde\Delta(1-e))=0$
for all $\omega\in A^*$ with $\omega\geq 0$.
This implies that $(\widetilde\Delta(1-e))(1\otimes e)=0$.
Because~$e$ is central, we also get $(\widetilde\Delta(x^*(1-e)x))(1\otimes e)=0$ for all $x\in A^{**}$ satisfying
$\widetilde\psi(x^*x)<\infty$.
If we apply $\widetilde\psi$, we get by using right invariance, that $\widetilde\psi(x^*(1-e)x)e=0$.
So, $\widetilde\psi(x^*(1-e)x)=0$ because~$e$ is non-zero.
Then we get $f(1-e)=0$.
A~similar argument will give $e(1-f)=0$ and therefore $e=f$.
\end{proof}

Compare this result with Theorem~3.8 in~\cite{K-V2}.

It follows immediately from this result that all invariant weights have the same support.
This implies that we have a~single von Neumann algebra.
Indeed, we can consider the associated von Neumann algebra~$M$, def\/ined as $A^{**}e$, where~$e$ is the support of the
non-trivial, densely def\/ined, lower semi-continuous, central invariant weights.

All the previous results lead to the following which is the main theorem of this appendix.

\begin{Theorem}\label{A.8}
Let~$A$ be a~$C^*$-algebra with a~comultiplication~$\Delta$.
Assume that there exist faithful, densely defined lower semi-continuous weights~$\varphi$ and~$\psi$, with central
support and resp.\ left and right invariant.
Let $M=A^{**}e$ where~$e$ is the support of these weights.
Consider the extension $\widetilde\Delta$ of~$\Delta$ to $A^{**}$ $($as in the proof of Proposition~{\rm \ref{A.6})}.
Then restrict to~$M$ and define $\Delta_0(x)=\widetilde\Delta(x)(e\otimes e)$ for $x\in M$.
This is a~comultiplication on the von Neumann algebra~$M$.
The restrictions to~$M$ of the extensions $\widetilde\varphi$ and $\widetilde\psi$ $($as in Proposition~{\rm \ref{A.3})} are a~left
and a~right Haar weight on $(M,\Delta_0)$ $($as in Definition~{\rm \ref{1.4}} and further in Section~{\rm \ref{Section2})}.
So, the pair $(M,\Delta_0)$ is a~locally compact quantum group in the sense of Definition~{\rm \ref{2.1}}.
\end{Theorem}

Remark that $(\widetilde\Delta(1-e))(1\otimes e)=0$ as we saw in the proof of Proposition~\ref{A.7}.
Then $(\widetilde\Delta(e))(1\otimes e)=1\otimes e$ and so $\Delta_0(e)=e\otimes e$.
This guarantees that the comultiplication $\Delta_0$ on~$M$ is unital.
This is also needed to show that $\Delta_0$ is still coassociative.

So, roughly speaking, Theorem~\ref{A.8} says that starting with a~$C^*$-algebra with a~comultiplication and nice
invariant weights, we can associate a~locally compact quantum group (in the von Neumann algebraic sense).
If we combine Proposition~\ref{A.2} with Proposition~\ref{A.5}, we see that also conversely, given a~locally compact
quantum group, we can associate a~pair of a~$C^*$-algebra and a~comultiplication with nice invariant weights.

What happens when we perform these two operations, one after the other? Do we get back the original?

First, start with a~locally compact quantum group $(M,\Delta)$.
Consider the $C^*$-algebra~$A$ as in Proposition~\ref{A.2} and restrict~$\Delta$ as in Proposition~\ref{A.7}.
Then it is rather straightforward to show that the construction in Theorem~\ref{A.8} will yield the original pair $(M,\Delta)$.

On the other hand, take a~$C^*$-algebra~$A$ with a~comultiplication~$\Delta$ and nice invariant weights.
Consider the pair $(M,\Delta_0)$ as in Theorem~\ref{A.8}.
We have the following lemma.

\begin{Lemma}\label{A.9}
Let~$W$ be the left regular representation for the pair $(M,\Delta_0)$.
Then the norm closure of the space
\begin{gather*}
\big\{(\iota\otimes \omega)W \,|\, \omega\in \mathcal B(\mathcal H_\varphi)_*\big\}
\end{gather*}
is the same as the norm closure of the space
\begin{gather*}
\SP\big\{(\iota\otimes \omega)\Delta(a) \,|\, a\in A,\; \omega\in A^*\big\}.
\end{gather*}
\end{Lemma}

This result is essentially found along with the proof of Proposition~\ref{1.21}.

It is not clear whether or not, this space will be all of~$A$.
So, in order to recover the original $C^*$-algebra, we need the extra {\it density conditions} as they are found in the
original def\/inition of a~locally compact quantum group in the $C^*$-algebra setting (as given by Kustermans and Vaes in~\cite{K-V2}).
We recall the def\/inition here.

\begin{Definition}\label{A.10}
Let~$A$ be a~$C^*$-algebra and~$\Delta$ a~comultiplication on~$A$ (as in Def\/inition~\ref{A.1}).
Assume that the spaces
\begin{gather*}
 \SP\{(\iota\otimes \omega)\Delta(a) \,|\, a\in A,\; \omega\in A^*\},
\qquad
 \SP\{(\omega\otimes\iota)\Delta(a) \,|\, a\in A,\; \omega\in A^*\}
\end{gather*}
are (norm) dense in~$A$.
Assume that there exist faithful, densely def\/ined lower semi-continuous weights~$\varphi$ and~$\psi$ on~$A$, both with
central support and resp.\ left and right invariant.
Then $(A,\Delta)$ is called a~locally compact quantum group in the $C^*$-algebraic sense.
\end{Definition}

It is only for such a~pair $(A,\Delta)$ that we have a~complete equivalence of the $C^*$-algebraic and von Neumann
algebraic setting.
If the density conditions are not fulf\/illed, in some sense, the $C^*$-algebra might be too big.
As we see from the previous discussion, we can pass to a~smaller $C^*$-algebra, invariant under the comultiplication,
satisfying the conditions of Def\/inition~\ref{A.10}.
However, it should be remarked that at present, there are no (obvious) examples of this phenomenon.

Next, one needs to prove properties for a~locally compact quantum group $(A,\Delta)$ with a~$C^*$-algebra~$A$ and
a~comultiplication~$\Delta$ as in Def\/inition~\ref{A.10}.
These properties will follow from the ones proven for a~locally compact quantum group $(M,\Delta)$ as def\/ined in
Def\/inition~\ref{2.1}.
We have already given some indications (about the stronger density conditions, the invariance of the $C^*$-algebra under
the modular automorphisms, the scaling automorphisms, the unitary antipode, \dots).
But more results need to be considered.
One can e.g.\ show relatively easy that $\delta^{it}\in M(A)$ for all~$t$ where~$\delta$ is the modular element from Theorem~\ref{2.11}.
There is also the construction of the dual $\widehat A$.
It can be obtained either by applying the procedure in this appendix and f\/ind $\widehat A$ from $\widehat M$ (as in Propositions~\ref{A.2}
and~\ref{A.5})
or directly as the norm closure of the space $\left\{(\omega\otimes\iota)W \,|\, \omega\in M_* \right\}$ (as
indicated in the~remark following Proposition~\ref{3.3}~-- compare also with the result in Lemma~\ref{A.9}).

So far about the comparison of our approach in this paper, with the $C^*$-algebraic approach by Kustermans and Vaes.
Let us now also compare (brief\/ly) with the approach of Masuda, Nakagami and Woronowicz.

There is f\/irst the original work by Masuda and Nakagami~\cite{M-N} where locally compact quantum groups are studied in
the von Neumann algebra framework.
Then, there there is the more recent work by these two authors and Woronowicz~\cite{M-N-W} where locally compact quantum
groups are studied in the $C^*$-algebra setting.
We will not say anything more about the f\/irst paper as in some sense, the second one can be seen as an improvement of
the f\/irst one.
We refer to the introduction of~\cite{M-N-W} for a~comparison of the two papers.

The main dif\/ference between the set of axioms in~\cite{M-N-W} and those formulated by Kustermans and Vaes is that in the
f\/irst case, the antipode and its polar decomposition are assumed whereas in~\cite{K-V2}, the antipode is constructed and
its properties are proven.
The same holds for our approach.

On the other hand, for this approach, as we have seen, we need to assume the existence of both a~left and a~right Haar weight.
This is not the case in the approach of~\cite{M-N-W} where only a~left Haar weight~$\varphi$ is assumed.
However, the unitary antipode~$R$ (see Def\/inition~1.5 in~\cite{M-N-W}) is part of the axioms and from all the axioms, it
is not so hard to obtain that~$R$ f\/lips the coproduct (see Proposition~2.6 in~\cite{M-N-W}).
Then the composition $\varphi\circ R$ gives a~right Haar weight.

Moreover, in examples, it is often rather easy to see what this unitary antipode should be and to verify that it f\/lips the coproduct.
Hence, usually, only one Haar weight is constructed explicitly while the other one is obtained by composing it with this
candidate for the unitary antipode.
Therefore it is in general more easy to verify the axioms of Kustermans and Vaes.

Another (minor) dif\/ference is that in~\cite{K-V2} weaker density conditions are needed (see earlier).
Also a~stronger invariance condition is proven by Kustermans and Vaes.

It is well-known that Masuda, Nakagami and Woronowicz started with their work on locally compact quantum groups, earlier
than Kustermans and Vaes but that it took many years before their results were published.
In my opinion, it is clear that the approach of Kustermans and Vaes, and also our approach in this paper, is better than their approach.
The axioms are more complicated and this does not really help to get the main results in a~quicker way.
On the other hand, the contribution of Masuda, Nakagami and Woronowicz is of great importance to the theory and their
paper contains interesting material and nice independent results.

Let us f\/inish this discussion by pointing out that it was in fact Kirchberg, at a~conference in Copenhagen, in~1992, who
f\/irst formulated the idea of generalizing the axioms of Kac algebras and thereby considering the polar decomposition of the antipode~\cite{K}.
As far as I know, his work was never published.

\subsection*{Acknowledgements}

Part of this work was done while I was on sabbatical in Trondheim (2002--2003).
First of all, I~would like to thank my colleague in Leuven, J.~Quaegebeur, who did some of my teaching in Leuven and
thus made it easier for me to go on sabbatical leave.
It is also a~pleasure to thank my colleagues of the University in Trondheim, M.~Landstad and C.~Skau, for their
hospitality during my stay and for giving me the opportunity to talk about my work in their seminar.
Part of this work was also done while visiting the University of Fukuoka (in~2004) and while visiting the University of
Urbana-Champaign (in~2005).
I~would also like to thank A.~Inoue, H.~Kurose and Z.-J.~Ruan for the hospitality during these (and earlier) visit(s).
I am also grateful to my coworkers in Leuven, especially J.~Kustermans and S.~Vaes, who have always been willing to
share their competence in this f\/ield.
This has been of great importance for my deeper understanding of the theory of locally compact quantum groups.
Again, without this, the present paper would not have been written.
Finally, I also wish to thank the organizers of the Special Program at the Fields Institute in June~2013 for giving me
the opportunity to present this work in a~series of lectures.

\pdfbookmark[1]{References}{ref}
\LastPageEnding

\end{document}